# EFFICIENCY IMPROVEMENTS IN INFERENCE ON STATIONARY AND NONSTATIONARY FRACTIONAL TIME SERIES[1]

### By P. M. Robinson

#### *London School of Economics*

We consider a time series model involving a fractional stochastic component, whose integration order can lie in the stationary/invertible or nonstationary regions and be unknown, and an additive deterministic component consisting of a generalized polynomial. The model can thus incorporate competing descriptions of trending behavior. The stationary input to the stochastic component has parametric autocorrelation, but innovation with distribution of unknown form. The model is thus semiparametric, and we develop estimates of the parametric component which are asymptotically normal and achieve an $M$-estimation efficiency bound, equal to that found in work using an adaptive LAM/LAN approach. A major technical feature which we treat is the effect of truncating the autoregressive representation in order to form innovation proxies. This is relevant also when the innovation density is parameterized, and we provide a result for that case also. Our semiparametric estimates employ nonparametric series estimation, which avoids some complications and conditions in kernel approaches featured in much work on adaptive estimation of time series models; our work thus also contributes to methods and theory for nonfractional time series models, such as autoregressive moving averages. A Monte Carlo study of finite sample performance of the semiparametric estimates is included.

**1. Introduction.** This paper obtains efficient parameter estimates in stationary or nonstationary, possibly fractional, time series. Consider a regression model given by

$$(1.1) \qquad y_t = \mu^T z_t + x_t, \qquad t \in \mathbb{Z},$$

Received June 2003; revised September 2004.

[1]Supported by a Leverhulme Trust Personal Research Professorship and ESRC Grant R000239936.

where $\mathbb{Z} = \{t : t = 0, \pm 1, \dots\}$, $z_t$ is a deterministic $q \times 1$ vector sequence, $\mu$ is an unknown $q \times 1$ vector, $T$ denotes transposition, $x_t$ is a zero-mean stochastic process and $y_t$ is an observable sequence. Any nonstationarity in the mean of $y_t$ would be due to $z_t$, nonstationarity in variance to $x_t$, but cases when $\mu^T z_t$ is a priori constant and $x_t$ is stationary are also of interest.

To describe $x_t$, denote by $B$ the back-shift operator, so $Bx_t = x_{t-1}$, and denote by $\Delta = 1 - B$ the difference operator; formally, for all real $d$

$$\Delta^{-d} = \sum_{j=0}^{\infty} \Delta_j(d) B^j, \qquad \Delta_j(d) = \frac{\Gamma(j+d)}{\Gamma(d)\Gamma(j+1)},$$

with $\Gamma$ denoting the gamma function such that $\Gamma(d) = \infty$ for $d = 0, -1, -2, \dots$, and $\Gamma(0)/\Gamma(0) = 1$. Assume the sequence $x_t$ is given by

$$(1.2) \qquad x_t = \Delta^{-m_0} v_t^{\#}, \qquad t \in \mathbb{Z},$$

where $m_0$ is a nonnegative integer,

$$(1.3) \qquad v_t^{\#} = v_t \mathbb{1}(t \geq 1), \qquad t \in \mathbb{Z},$$

for $\mathbb{1}(\cdot)$ the indicator function, and

$$(1.4) \qquad v_t = \Delta^{-\zeta_0} u_t, \qquad t \in \mathbb{Z},$$

for $|\zeta_0| < \frac{1}{2}$, with $u_t$ a zero-mean covariance stationary process with absolutely continuous spectral distribution function and spectral density $f(\lambda)$ that is at least positive and finite for all $\lambda$.

The process $v_t$ is then also covariance stationary, having "long memory" for $\zeta_0 > 0$, "short memory" for $\zeta_0 = 0$ and "negative memory" for $\zeta_0 < 0$. When $m_0 = 0$, we have $x_t = v_t^{\#} = v_t$ for $t \geq 1$. When $m_0 \geq 1$, $x_t$ "integrates" $v_t^{\#}$, and the truncation in (1.2) implies that $x_t$ has variance that is finite, albeit evolving with $t$. Putting $\xi_0 = m_0 + \zeta_0$, $x_t$ is well defined for

$$(1.5) \qquad \xi_0 \in S \subset \{\xi : -\tfrac{1}{2} < \xi < \infty, \xi \neq \tfrac{1}{2}, \tfrac{3}{2}, \dots\}.$$

The requirement $\xi_0 > -\frac{1}{2}$ excludes noninvertible processes, and the final qualification in (1.5) excludes $\xi_0$ that cannot be reduced to the stationary/invertible region $(-\frac{1}{2}, \frac{1}{2})$ by integer differencing. Alternative definitions of nonstationary fractional $x_t$ are available, for example, $\Delta^{-\xi_0} u_t^{\#}$.

Suppose $\xi_0$ is unknown; $m_0$ may also be unknown. Suppose $u_t$ is assumed to have parametric autocorrelation,

$$(1.6) \qquad f(\lambda) = \frac{\sigma_0^2}{2\pi} |\beta(e^{i\lambda}; \nu_0)|^2, \qquad \lambda \in (-\pi, \pi],$$

such that $\operatorname{cov}(u_0, u_j) = \int_{-\pi}^{\pi} f(\lambda) \cos(j\lambda) \, d\lambda$, $j \in \mathbb{Z}$, $\beta(s; \nu)$ is a smooth given function of complex-valued $s$ and the column-vector $\nu \in V \subset \mathbb{R}^{p_1-1}$, $p_1 \geq 1$, satisfying

$$(1.7) \qquad \beta_0(\nu) = 1, \qquad \beta(s; \nu) \neq 0, \qquad |s| \leq 1, \nu \in V,$$

where $\beta_j(\nu) = \int_{-\pi}^{\pi} \beta(e^{i\lambda}; \nu) \cos(j\lambda)\,d\lambda$, and $\nu_0 \in V$ and $\sigma_0^2 > 0$ are unknown. Then $\sigma_0^2$ is the variance of the one-step-ahead prediction error of the best linear predictor for $u_t$. For example, $u_t$ can be a standardly parameterized autoregressive moving average (ARMA) process of autoregressive (AR) order $p_{11}$ and moving average (MA) order $p_{12}$, such that $p_1 - 1 \le p_{11} + p_{12} < \infty$; when $\nu_0$ consists precisely of the AR and MA coefficients we have $p_{11} + p_{12} = p_1 - 1$; otherwise the coefficients obey prior restrictions. We call $v_t$ a FARIMA$(p_{11}, \zeta_0, p_{12})$, and $x_t$ a FARIMA$(p_{11}, \xi_0, p_{12})$. Whereas $v_t$ is stationary, due to the truncation (1.2) $x_t$ is nonstationary even when $\xi_0 < \frac{1}{2}$ (it could be called "asymptotically stationary" then). The case when $x_t = v_t$ for all $t \in \mathbb{Z}$, so $x_t$ is stationary, can be dealt with similarly, but we impose the truncation in (1.2) for all $m_0 \ge 0$ for the sake of a unified presentation. The set $V$ is contained in the "stationary and invertible region." The case $p_1 = 1$ means $\nu_0$ is empty, and if $\beta \equiv 1$, $x_t$ is a FARIMA$(0, \xi_0, 0)$. An alternative model for $u_t$ is due to Bloomfield [4].

The main focus of the paper is estimation of $\theta_{01} = (\xi_0, \nu_0^T)^T$, and we restrict to a specialized form of $z_t$ in (1.1):

$$(1.8) \qquad z_t = (t^{\tau_1}, \dots, t^{\tau_q})^T \mathbb{1}(t \ge 1), \qquad \tau_1 < \tau_2 < \cdots < \tau_q,$$

where the $\tau_j$ are real valued. Debate has centered on the origin—deterministic or stochastic—of nonstationarity in time series. A notable feature is competition at low frequencies, and given the fractional model for $x_t$ this is most neatly expressed by (1.8). Some components of $z_t$ may have negligible effect on fractionally differenced $y_t$. Denote by $\mu_j$ the $j$th element of $\mu$ and $\mathcal{T}_1 = \{j : \tau_j < \xi_0 - \frac{1}{2}\}$, $\mathcal{T}_2 = \{j : \tau_j = \xi_0\}$, $\mathcal{T}_3 = \{j : \xi_0 - \frac{1}{2} \le \tau_j < \xi_0; \tau_j > \xi_0\}$, where any of these sets can be empty. We cannot estimate $\mu_j$ for $j \in \mathcal{T}_1$, and do not discuss estimation of $\mu_j$ for $j \in \mathcal{T}_2$. Write $s_t = \sum_{j \in \mathcal{T}_1} \mu_j t^{\tau_j}$ and for $p_2 = \# \mathcal{T}_3 \le q$ introduce the $p_2 \times 1$ vectors $z_{2t}$ and $\theta_{02}$, whose $j$th elements are the elements of $z_t$ and $\mu$ whose index is the $j$th largest element of $\mathcal{T}_3$. It will be convenient to write $z_{2t} = (t^{\chi_1}, \dots, t^{\chi_{p_2}})^T$, where the $\chi_j$ are appropriate elements $\tau_j$, and satisfy $\frac{1}{2} \le \chi_1 < \cdots < \chi_{p_2}$. We can write (1.1) as

$$(1.9) \qquad y_t = s_t + \mu^* t^{\xi_0} + \theta_{02}^T z_{2t} + x_t,$$

where $\mu^* = 0$ if $\tau_j \ne \xi_0$ for all $j$.

We discuss estimation of $\theta_{02}$, along with $\theta_{01}$. For this we require that the $\tau_j$, $j \in \mathcal{T}_3$, are known. The boundary case of $\mathcal{T}_3$, $\tau_j = \xi_0 - \frac{1}{2}$, thus strictly implies $\xi_0$ is known, but this provision is instead designed to cover a situation in which $\tau_j < \xi_0 - \frac{1}{2}$ for all $j \in \mathcal{T}_1$ is anticipated, with $\xi_0$ unknown, but in fact $\tau_j = \xi_0 - \frac{1}{2}$ for some $j$. For $\theta_1 = (\xi, \nu^T)^T \in S \times V$, introduce the function $\alpha(s; \theta_1) \colon \mathbb{R} \times \mathbb{R}^{p_1} \to \mathbb{R}$, and consider $\alpha(s; \theta_1^{(-)})$, where $\theta_1^{(-)} = (0, \nu^T)^T$, such that

$$(1.10) \qquad \alpha(s; \theta_1) = (1 - s)^\xi \alpha(s; \theta_1^{(-)}).$$

Take $\alpha(s; \theta_1^{(-)}) = \beta(s; \nu)^{-1}$ for $|s| \le 1$, $\nu \in V$, and note that $\int_{-\pi}^{\pi} \alpha(e^{i\lambda}; \theta_1^{(-)}) \, d\lambda = 1$, $\nu \in V$. From (1.6) and (1.7), $u_t$ has one-sided AR representation

$$\alpha(B; \theta_{01}^{(-)}) u_t = \sigma_0 \varepsilon_t, \qquad t \in \mathbb{Z}, \tag{1.11}$$

where $\theta_{01}^{(-)} = (0, \nu_0^T)^T$, and the $\varepsilon_t$ are uncorrelated with zero mean and unit variance. Introduce square-summable coefficients $\alpha_j(\theta_1)$ in the expansion

$$\alpha(s; \theta_1) = \sum_{j=0}^{\infty} \alpha_j(\theta_1) s^j, \qquad |s| \le 1, \xi \in S, \nu \in V, \tag{1.12}$$

so $\alpha_0(\theta_1) \equiv 1$. For given $\theta = (\theta_1^T, \theta_2^T)^T$, define the computable

$$e_t(\theta) = \sum_{j=0}^{t-1} \alpha_j(\theta_1)(y_{t-j} - \theta_2^T z_{2, t-j}),$$

$$E_t(\theta) = e_t(\theta) - \frac{1}{n} \sum_{t=1}^{n} e_t(\theta), \qquad t \ge 1, \tag{1.13}$$

the latter being proxies for $\sigma_0 \varepsilon_t$, with $s_t$ ignored in $e_t(\theta)$ because it is anticipated to have negligible effect, and $\mu^* t^{\xi_0}$ ignored in view of the mean-correction in $E_t(\theta)$.

Given observations $y_t$, $t = 1, \ldots, n$, define

$$Q_\rho(\theta, \theta_3) = \frac{1}{n} \sum_{t=1}^{n} \rho(E_t(\theta)/\tilde{\sigma}; \theta_3), \tag{1.14}$$

for an $n^{1/2}$-consistent estimate $\tilde{\sigma}$ of $\sigma_0$, a given nonnegative function $\rho : \mathbb{R} \times \mathbb{R}^{p_3} \Rightarrow \mathbb{R}$ and any admissible value $\theta_3$ of an unknown $p_3 \times 1$ parameter vector $\theta_{03}$; $\theta_3$ may be empty, as when $\rho(s; \theta_3) = s^2$. Consider the estimate $(\bar{\theta}_\rho^T, \bar{\theta}_{3\rho}^T) = \arg \min_{\Theta \times \Theta_3} Q_\rho(\theta, \theta_3)$, for compact sets $\Theta \in \mathbb{R}^p$, $\Theta_3 \in \mathbb{R}^{p_3}$. One anticipates (see, e.g., Martin's [24] discussion of $M$-estimates of ARMA models) that under suitable conditions $\bar{\theta}_\rho, \bar{\theta}_{3\rho}$ are asymptotically independent and the asymptotic variance matrix of $\bar{\theta}_\rho$ depends on $\rho$ only through the scalar factor $\mathcal{H} = \int \rho'(s)^2 g(s) \, ds / \{\int \rho''(s) g(s) \, ds\}^2$, where the prime indicates differentiation, double-prime indicates twice differentiation and reference to $\theta_{03}$ is suppressed. If integration by parts can be conducted, this and the Schwarz inequality indicate that $\mathcal{H} \ge \mathcal{J}^{-1}$, defining the information

$$\mathcal{J} = \int \psi(s)^2 g(s) \, ds \tag{1.15}$$

and the score function

$$\psi(s) = -g'(s)/g(s). \tag{1.16}$$

The lower bound is attained by $\bar{\bar{\theta}}_{\log g}$, and the paper obtains estimates that are efficient in the sense of having the same asymptotic variance as $\bar{\bar{\theta}}_{\log \rho}$. In Theorem 2 of Section 3 we justify such an estimate on the basis of known $g(s; \theta_3)$. If $g$ is misspecified, not only will the estimate not be efficient but it may even be inconsistent. Our main result is Theorem 1 of Section 3, which justifies efficient semiparametric estimates, in which the density of $\varepsilon_t$ is nonparametric. These estimates are adaptive in the sense of Stone [28] and are described in the following section. Section 4 describes a Monte Carlo study of finite sample behavior of the semiparametric estimates. Section 5 attempts to place the work in perspective, relative to the literature. Section 6 presents the main proof details, which use a series of lemmas that make up Section 7. Some of these, such as Lemmas 1, 2, 7, 8, 13, 15 and 16, may be useful in other work. A principal technical feature is our handling of the approximation of the $\sigma_0 \varepsilon_t$ in (1.11) by the $e_t(\theta_0)$ defined by (1.13), a delicate matter in fractional models.

**2. Semiparametric estimates.** As in much adaptive estimation literature we take an approximate Newton step from an initial consistent estimate $\tilde{\theta}$ of $\theta_0$, with the same rate of convergence as $\bar{\bar{\theta}}_{\log g}$. This requires estimating $\psi(s)$. We employ an approach developed by Beran [2] and Newey [25]. Beran [2] proposed a series estimate of $\psi(s)$ [with respect to innovations in an AR($p$) model] that employs integration by parts. His estimate of $\psi(s)$ was actually not a smoothed nonparametric one because he fixed the number of terms, $L$, in the series. Newey [25] allowed $L$ to increase slowly with $n$, in adapting to error distribution of unknown form in cross-sectional regression.

Let $\phi_\ell(s)$, $\ell = 1, 2, \ldots$, be a sequence of given, continuously differentiable functions. For $L \geq 1$, scalar $h_t$, $t = 1, \ldots, n$, and $h = (h_1, \ldots, h_n)^T$, define $\phi^{(L)}(h_t) = (\phi_1(h_t), \ldots, \phi_L(h_t))^T$, $\Phi^{(L)}(h_t) = \phi^{(L)}(h_t) - n^{-1} \sum_{s=1}^n \phi^{(L)}(h_s)$, $\phi'^{(L)}(h_t) = (\phi_1'(h_t), \ldots, \phi_L'(h_t))^T$ and

$$W^{(L)}(h) = n^{-1} \sum_{t=1}^n \Phi^{(L)}(h_t) \Phi^{(L)}(h_t)^T,$$

$$w^{(L)}(h) = n^{-1} \sum_{t=1}^n \phi'^{(L)}(h_t),$$

$$\hat{a}^{(L)}(h) = W^{(L)}(h)^{-1} w^{(L)}(h),$$

$$\psi^{(L)}(h_t; \hat{a}^{(L)}(h)) = \hat{a}^{(L)}(h)^T \Phi^{(L)}(h_t).$$

With $E(\theta) = (E_1(\theta), \ldots, E_n(\theta))^T$ define

$$\tilde{\psi}_t^{(L)}(\theta, \sigma) = \psi^{(L)}(E_t(\theta)/\sigma; \hat{a}^{(L)}(E(\theta)/\sigma)),$$

where it will follow from our conditions that in a neighborhood of $\theta_0, \sigma_0$, $W^{(L)}(E(\theta)/\sigma)$ is nonsingular with probability approaching 1 as $n \to \infty$.

We then compute the $\tilde{\psi}_t^{(L)}(\tilde{\theta}, \tilde{\sigma})$. Following Beran [2] and Newey [25] we have approximated $\psi(\varepsilon_t)$ by $\sum_{\ell=1}^L a_\ell \{\phi_\ell(\varepsilon_t) - E\phi_\ell(\varepsilon_t)\}$ [imposing the restriction $E\psi(\varepsilon_t) = 0$], noted that (under conditions to be given) integration by parts implies $E\{\phi^{(L)}(\varepsilon_t)\psi(\varepsilon_0)\} = E\{\phi^{(L)}(\varepsilon_t)\}$, estimated $(a_1, \ldots, a_L)^T$ by $a^{(L)}(E(\tilde{\theta})/\tilde{\sigma})$, and then $\psi(\varepsilon_t)$ by $\tilde{\psi}_t^{(L)}(\tilde{\theta}, \tilde{\sigma})$.

Define [see (1.10)–(1.13)]

$$e_t'(\theta) = (\partial/\partial\theta)e_t(\theta) = (e_{t1}'(\theta)^T, e_{t2}'(\theta)^T)^T,$$

where

$$e_{t1}'(\theta) = \alpha'(B; \theta_1)(y_t - \theta_2^T z_{2t}), \qquad e_{t2}'(\theta) = -\alpha(B; \theta_1)z_{2t},$$

with

$$(2.1) \qquad \begin{aligned} \alpha'(s; \theta_1) &= (\partial/\partial\theta_1)\alpha(s; \theta_1) = (1 - s)^\xi \alpha(s; \theta_1^{(-)})\gamma(s; \nu), \\ \gamma(s; \nu) &= [\log(1 - s), \{(\partial/\partial\nu)^T\alpha(s; \theta_1^{(-)})\}/\alpha(s; \theta_1^{(-)})]^T. \end{aligned}$$

Define

$$E_{ti}'(\theta) = e_{ti}'(\theta) - n^{-1}\sum_{s=1}^n e_{si}'(\theta), \qquad i = 1, 2,$$

$$r_{Li}(\theta, \sigma) = \sum_{t=1}^n \tilde{\psi}_t^{(L)}(\theta, \sigma)E_{ti}'(\theta), \qquad R_i(\theta) = \sum_{t=1}^n E_{ti}'(\theta)E_{ti}'(\theta)^T, \qquad i = 1, 2,$$

$$\mathcal{J}_L(\theta, \sigma) = n^{-1}\sum_{t=1}^n \tilde{\psi}_t^{(L)}(\theta, \sigma)^2.$$

Estimate $\theta_{01}, \theta_{02}$ by

$$(2.2) \qquad \hat{\theta}_i = \tilde{\theta}_i + \{R_i(\tilde{\theta})\mathcal{J}_L(\tilde{\theta}, \tilde{\sigma})\}^{-1}r_{Li}(\tilde{\theta}, \tilde{\sigma}), \qquad i = 1, 2,$$

respectively, for $\tilde{\theta} = (\tilde{\theta}_1^T, \tilde{\theta}_2^T)^T$.

As in [25] we restrict to $\phi_\ell(s)$ satisfying

$$(2.3) \qquad \phi_\ell(s) = \phi(s)^\ell,$$

for a smooth function $\phi(s)$. Examples are

$$(2.4) \qquad \phi(s) = s,$$

$$(2.5) \qquad \phi(s) = s(1 + s^2)^{-1/2}.$$

Our conditions require $L$ to increase very slowly with $n$, and allow the increase to be arbitrarily slow; in practice, for moderate $n$, (2.2) might be computed for a few small integers $L$, starting with $L = 1$. Recursive formulas are available, using partitioned regression, such that the elements of $W^{(L)}(E(\tilde{\theta})/\tilde{\sigma})$, $w^{(L)}(E(\tilde{\theta})/\tilde{\sigma})$ can be used in computing $\tilde{\psi}_t^{(L+1)}(\tilde{\theta}, \tilde{\sigma})$.

**3. Main results.** We introduce the following regularity conditions. Throughout the paper $C$ denotes a finite but arbitrarily large constant.

ASSUMPTION A1. The sequence $y_t$ is generated by (1.1) with $x_t$ generated by (1.2)–(1.4) and (1.11), where the $\varepsilon_t$ are independent and identically distributed (i.i.d.) with zero mean and variance 1, and $z_t$ is given by (1.8).

ASSUMPTION A2. Either:

(a) $E\varepsilon_0^4 < \infty$; or
(b) for some $\omega > 0$ the moment generating function $E(e^{t|\varepsilon_0|^\omega})$ exists for some $t > 0$; or
(c) $\varepsilon_0$ is almost surely bounded.

ASSUMPTION A3. $\varepsilon_0$ has density $g(s)$ that is differentiable and

$$0 < \mathcal{J} < \infty,$$

where $\mathcal{J}$ is defined in (1.15).

ASSUMPTION A4. The sentence including (1.6) and (1.7) is true, $\nu_0$ is an interior point of $V$, and in a neighborhood $\mathcal{N}$ of $\nu_0$, $\alpha(s; \theta_1^{(-)}) = \beta(s; \nu)^{-1}$ is thrice continuously differentiable in $\nu$ for $|s| = 1$ and

$$\sum_{j=1}^{\infty} j^3 \left\{ |\beta_j(\nu_0)| + \sup_{\mathcal{N}} |\alpha_j(\theta_1^{(-)})| + \sup_{\mathcal{N}} |\alpha_j^{(k)}(\theta_1^{(-)})| \right. $$
$$\left. + \sup_{\mathcal{N}} |\alpha_j^{(k,\ell)}(\theta_1^{(-)})| + \sup_{\mathcal{N}} |\alpha_j^{(k,\ell,m)}(\theta_1^{(-)})| \right\} < \infty,$$

for all $k, \ell, m = 1, \ldots, p_1 - 1$, where $\alpha_j(\theta_1^{(-)})$ is defined by (1.10), (1.12) and $\alpha_j^{(k)}(\theta_1^{(-)}) = (\partial/\partial\nu_k)\alpha_j(\theta_1^{(-)})$, $\alpha_j^{(k,\ell)}(\theta_1^{(-)}) = (\partial/\partial\nu_\ell)\alpha_j^{(k)}(\theta_1^{(-)})$, $\alpha_j^{(k,\ell,m)}(\theta_1^{(-)}) = (\partial/\partial\nu_m)\alpha_j^{(k,\ell)}(\theta_1^{(-)})$, $\nu_k$ being the $k$th element of $\nu$.

ASSUMPTION A5. For all $(p_1 - 1) \times 1$ nonnull vectors $\lambda$, $\lambda^T \{(\partial/\partial\nu)\alpha(e^{i\lambda}; \theta_{01}^{(-)})\}\beta(e^{i\lambda}; \nu_0) \neq 0$ on a subset of $(-\pi, \pi]$ of positive measure.

ASSUMPTION A6.

$$0 < \sigma_0^2 < \infty.$$

ASSUMPTION A7.
$$n^{1/2}(\tilde{\theta}_1 - \theta_{01}) = O_p(1), \qquad D_n(\tilde{\theta}_2 - \theta_{02}) = O_p(1), \qquad n^{1/2}(\tilde{\sigma}^2 - \sigma_0^2) = O_p(1),$$

where

$$D_n = \text{diag}\{n^{\chi_1-\xi_0+1/2}\mathbb{1}(\chi_1-\xi_0 > -\tfrac{1}{2}) + (\log n)^{1/2}\mathbb{1}(\chi_1-\xi_0 = -\tfrac{1}{2}),$$
$$n^{\chi_2-\xi_0+1/2}, \ldots, n^{\chi_{p_2}-\xi_0+1/2}\}.$$

ASSUMPTION A8.   $\phi_\ell(s)$ satisfies (2.3), where $\phi(s)$ is strictly increasing and thrice continuously differentiable and is such that, for some $\kappa \geq 0$, $K < \infty$,

$$(3.1) \qquad\qquad |\phi(s)| \leq \mathbb{1}(|s| \leq 1) + |s|^\kappa \mathbb{1}(|s| > 1),$$

$$(3.2) \qquad |\phi'(s)| + |\phi''(s)| + |\phi'''(s)| \leq C(1 + |\phi(s)|^K).$$

ASSUMPTION A9.

$$(3.3) \qquad\qquad L \to \infty \qquad \text{as } n \to \infty$$

and either:

(a)

$$(3.4) \quad \liminf_{n\to\infty}\left(\frac{\log n}{L}\right) > 8\{\log\eta + \max(\log\varphi, 0)\} \simeq 7.05 + 8\max(\log\varphi, 0);$$

or

(b)

$$(3.5) \qquad\qquad \liminf_{n\to\infty}\left(\frac{\log n}{L\log L}\right) > \max\left(\frac{8\kappa}{\omega}, \frac{4\kappa(\omega+1)}{\omega}\right);$$

or

(c)

$$(3.6) \qquad\qquad \liminf_{n\to\infty}\left(\frac{\log n}{L\log L}\right) > 4\kappa,$$

where

$$\eta = 1 + 2^{1/2} \simeq 2.414$$

and

$$\varphi = \frac{1 + |\phi(s_1)|}{\phi(s_2) - \phi(s_1)},$$

$[s_1, s_2]$ being an interval on which $g(s)$ is bounded away from zero.

REMARK 1.   Parts (a), (b) and (c) of Assumption A2 increase in strength and entail trade-offs with Assumptions A8 and A9. When $\kappa = 0$ in Assumption A8, so $\phi(s)$ is bounded, (a) of Assumption A2 and (a) of Assumption A9 suffice; a finite fourth moment seems hard to avoid in dealing with the deviations $e_t(\theta_0) - \sigma_0 \varepsilon_t$. Part (b) of Assumption A2 holds with $\omega = 1$ for Laplace $\varepsilon_t$ and with $\omega = 2$ for Gaussian $\varepsilon_t$. We require (b) of Assumption A2 when $\kappa > 0$ in Assumption A8, so $\phi(s)$ can be unbounded, and also (b) of Assumption A9. If (c) of Assumption A2 holds, then a fortiori we can have $\kappa > 0$ in Assumption A8, and can relax (b) of Assumption A9 to (c).

REMARK 2. Assumption A3 is virtually necessary.

REMARK 3. Assumption A4 is stronger than necessary, but is chosen for brevity of presentation and because it is readily checked for short memory and invertible AR ($\alpha$) and MA ($\beta$) filters arising in models of most practical interest, such as ARMA and Bloomfield [4] models, and in any case conditions on the short-memory component are of only secondary interest here. A property useful in several places (see in particular Lemma 13 of Section 7) that is ensured by Assumption A4 is as follows. A (possibly vector) sequence $\alpha_j$, $j \geq 0$, has property $P_r(d)$, $r \geq 0$, if

$$\|\alpha_j\| \leq C\{\log(j+2)\}^r (j+1)^{d-1},$$

$$\|\alpha_j - \alpha_{j+1}\| \leq C\{\log(j+2)\}^r (j+1)^{d-2}, \qquad j \geq 0,$$

where $\|\cdot\|$ denotes Euclidean norm. For $|s| \leq 1$ and $\theta_1^{(+)} = (\zeta, \nu^T)^T$, define square-summable $\pi_j(\theta_1^{(+)})$ such that

$$\pi(s; \theta_1^{(+)}) = (1-s)^{-\zeta} \beta(s; \nu) = \sum_{j=0}^{\infty} \pi_j(\theta_1^{(+)}) s^j, \qquad |\zeta| < \tfrac{1}{2}, \nu \in V.$$

Then, with $\theta_{01}^+ = (\zeta_0, \nu_0^T)^T$, $\pi_j(\theta_{01}^{(+)})$ has property $P_0(\zeta_0)$, $\alpha_j(\theta_{01}^{(+)})$ has property $P_0(-\zeta_0)$ and $(\partial/\partial\theta_1^{(+)T}) \alpha_j(\theta_{01}^{(+)})$ has property $P_1(-\zeta_0)$. This follows from Lemmas 11 and 12 of Section 7 on noting that, for $\alpha(s) = \sum_{j=0}^{\infty} \alpha_j s^j$, $\beta(s) = \sum_{j=0}^{\infty} \beta_j s^j$, the coefficient of $s^j$ in $\alpha(s)\beta(s)$ is $\sum_{k=0}^{j} \alpha_k \beta_{j-k}$, that the coefficients of $s^j$ in $(1-s)^{-d}$ and $-\log(1-s)$ are $\Delta_j(d)$ and $j^{-1}$, that $\pi(1; \theta_{01}^{(+)}) = 0$ for $\zeta_0 < 0$, and that $\alpha(1; \theta_{01}^{(+)}) = 0$, $(\partial/\partial\theta^{(+)T}) \alpha(1; \theta_{01}^{(+)}) = 0$ for $\zeta_0 > 0$.

REMARK 4. Assumption A5 is an identifiability condition, violated if, for example, $u_t$ is specified as an ARMA with both AR and MA orders overstated. Assumption A5, with Assumption A4, implies that

$$
\begin{aligned}
\Omega_1 &= \frac{1}{2\pi} \int_{-\pi}^{\pi} \gamma(e^{i\lambda}; \nu_0) \gamma(e^{-i\lambda}; \nu_0)^T \, d\lambda \\
&= \frac{1}{2\pi} \int_{-\pi}^{\pi} \left[ \begin{array}{c} \log|1 - e^{i\lambda}|^2 \\ 2\dfrac{\partial}{\partial\nu} \log|\beta(e^{i\lambda}; \nu_0)| \end{array} \right] \left[ \begin{array}{c} \log|1 - e^{i\lambda}|^2 \\ 2\dfrac{\partial}{\partial\nu} \log|\beta(e^{i\lambda}; \nu_0)| \end{array} \right]^T \, d\lambda
\end{aligned}
$$

(3.7)

is positive definite, with $\gamma$ given by (2.1). $\Omega_1$ is proportional to the inverse of the limiting covariance matrix of $\hat{\theta}_1$. We define also the corresponding matrix with respect to $\hat{\theta}_2$,

$$
\begin{aligned}
(3.8) \quad \Omega_2 = {} & \frac{\sigma_0^2}{2\pi}\beta(1;\nu_0)^2 \\
& \times \left( \frac{\{2(\chi_i-\xi_0)+1\}^{1/2}\{2(\chi_j-\xi_0)+1\}^{1/2}(\chi_i-\xi_0)(\chi_j-\xi_0)}{(\chi_i+\chi_j-2\xi_0+1)(\chi_i-\xi_0+1)(\chi_j-\xi_0+1)} \right),
\end{aligned}
$$

when $\chi_1 - \xi_0 > -\frac{1}{2}$, where the $(i,j)$th element of the matrix is displayed; because $((\chi_i+\chi_j-2\xi_0+1)^{-1})$ is a Cauchy matrix (see [17], page 30), and the inequalities in (1.8) hold, $\Omega_2$ is positive definite. The same is true when $\tau_j - \xi_0 = -\frac{1}{2}$ for some $j$, $\Omega_2$ being defined by replacing the $(1,1)$th element of the matrix in (3.8) by 1, and the other elements in the first row and column by zero.

REMARK 5. The middle part of Assumption A7 is likely to be satisfied by the least-squares estimate of $\theta_{02}$, under similar conditions to ours. A substantial literature justifies $\hat{\theta}_1$ satisfying Assumption A7; typically $\theta'_{02}z_{2t}$ is assumed constant a priori, but the results should go through more generally with $x_t$ replaced by least-squares residuals. Various estimates of $\theta_{01}$ (which we collectively call Whittle estimates) have been shown to be $n^{1/2}$-consistent and asymptotically $N(0,\Omega_1^{-1})$ when $0 \le \xi_0 < \frac{1}{2}$ under Gaussianity of $x_t$ (when they achieve the efficiency bound of Section 1 and are as good as maximum likelihood estimates), and under more general conditions (see, e.g., [6, 9, 11, 16]). The estimate minimizing (1.14) with $\rho(s) = s^2$ [usually with $E_t(\theta)$ replaced by $e_t(\theta)$] falls within this class. This estimate (used by Li and McLeod [21] for fractional models and Box and Jenkins [5] for ARMA ones) is sometimes called a conditional sum of squares (CSS) estimate (though it is based on formulas for the truncated AR representation rather than for the conditional expectation given the finite past record). Beran [1] argued that it has the same desirable asymptotic properties for $\xi_0 > \frac{1}{2}$, tying in with Robinson's [26] derivation of standard asymptotics for score tests, based on the same objective function, for unit root and more general nonstationary hypotheses against fractional alternatives. These authors employed a different definition of fractional nonstationarity from ours, but for our definition Velasco and Robinson [29] established the same properties for a Whittle estimate when $-\frac{1}{2} < \xi_0 < \frac{3}{4}$, and for a tapered version of this for $-\frac{1}{2} < \xi_0 < \infty$, though the tapering inflates asymptotic variance. They established consistency of their implicitly defined optimizer despite lack of uniform convergence over an admissible parameter set that includes a wide range of nonstationary values of $\xi$. Taking a Newton step from a previously established $n^{1/2}$-consistent estimate avoids repeating this kind of work. Velasco and Robinson's [29] estimate of $\sigma_0^2$ should satisfy the final part of Assumption A7 [with (a) sufficient within Assumption A2].

REMARK 6. When $\kappa = 0$ in Assumption A8, then $|\phi(s)| \le 1$ for all $s$, under (3.1); there would be no gain in generality by specifying $\phi$ to satisfy a larger finite bound. For $\kappa > 0$ we might take $\phi(s) = s^\kappa$; compare (2.4). The reason for imposing different bounds on $\phi(s)$ over $|s| \le 1$ and $|s| > 1$ is to allow possibly different rates of approach to zero and infinity. Assumption A8 is stronger than the corresponding assumption of Newey [25], and is driven by the presence of $e_t(\theta_0)$ for small $t$, when it does not approximate $\sigma_0 \varepsilon_t$; we prefer this to trimming out small $t$, which introduces further ambiguity. It is hard to think of reasons for choosing $\phi$ that do not satisfy (3.1), (3.2), which imply power-law bounds on $\phi'(s)$, $\phi''(s)$ and $\phi'''(s)$ as $s \to \infty$.

REMARK 7. The weakest of the conditions in Assumption A9, (a), can only apply when $\kappa = 0$ in Assumption A8, in which case $\log \varphi > 0$. Subject to this, the hope is that $s_1$ and $s_2$ exist such that $\varphi$ is arbitrarily close to 1, as when $g(s) > 0$ for all $s$; then the strict inequality in (3.4) applies with $\log \varphi = 0$. The mysterious constant $\eta$ is due to approximating $W^{(L)}$ in the proof in terms of the Cauchy matrix with $(i, j)$th element $\int_{-1}^{1} u^{i+j-2} du$ (see Lemma 7 of Section 7). Since $\phi$ is defined for negative and positive arguments, this seems more natural than Newey's [25] use of the Hilbert matrix ($\int_0^1 u^{i+j-2} du$) and affords some slight improvement over it due to the many zero elements in this Cauchy matrix; following a similar proof to that of Lemma 7 for the Hilbert matrix, $\eta$ would be replaced by $\eta^2 \simeq 5.828$. In fact, a constant such as $\eta$ does not arise in Newey's work because he is content with a slightly stronger condition than any in Assumption A9, $L \log L / \log n \to 0$, irrespective of whether or not $\phi$ is bounded, and without considering the impact of bounded $\varepsilon_t$. This is because he accepts a bound of form $L^{CL}$ at several points of his proof. Our slightly sharper bounds suggest that when $\phi$ is bounded it is effectively the denominator of $\psi^{(L)}$ (i.e., the inverse of $W^{(L)}$) that dominates, while when $\phi$ is unbounded the numerator dominates. In the former case, the slow $L$ corresponds to the notorious ill-conditioning of Cauchy–Hilbert matrices. One disadvantage of a bounded $\phi$ is that a larger $L$ might be needed to approximate an unbounded $\psi$, though our slightly milder condition on $L$ in Assumption A9(a) might help to justify this. Another is that it excludes (2.4), which "nests" the Gaussian case, though it would be possible to modify our theory to allow inclusion of $\phi_1(s) = s$, say, followed by polynomial $\phi_\ell$ (2.3) using bounded $\phi$ such as (2.5). Though the partly known nature of the bounds in Assumption A9 is interesting, and their reflection of other assumptions is intuitively reasonable in a relative sense, not only is the improvement over Newey's rate slight, but even after guessing $\omega$ and $\varphi$, no practical choices of $L$ in finite samples can be concluded; indeed the same asymptotic bounds result if any fixed integer is added to or subtracted from $L$. As in much other semiparametric work, no information toward an optimal choice of $L$ emerges; indeed, as in [25] there is no lower bound on $L$, and besides that it must increase with $n$.

THEOREM 1. *Let Assumptions* A1–A9 *hold, such that when $\kappa = 0$ Assumption* A2(a) *holds with Assumption* A9(a), *or when $\kappa > 0$ either Assumption* A2(b) *holds with Assumption* A9(b) *or Assumption* A2(c) *holds with Assumption* A9(c). *Then as $n \to \infty$, $n^{1/2}(\hat{\theta}_1 - \theta_{01})$ and $D_n(\hat{\theta}_2 - \theta_{02})$ converge in distribution to independent $N(0, \mathcal{J}^{-1}\Omega_1^{-1})$, $N(0, \mathcal{J}^{-1}\Omega_2^{-1})$ vectors, respectively, where the limiting covariance matrices are consistently estimated by $\{\mathcal{J}_L(\tilde{\theta}, \tilde{\theta})R_1(\tilde{\theta})/n\}^{-1}$, $\{\mathcal{J}_L(\tilde{\theta}, \tilde{\theta})D_n^{-1}R_2(\tilde{\theta})D_n^{-1}\}^{-1}$, respectively.*

To place Theorem 1 in perspective and to further balance the focus on Whittle estimation in the long-memory literature, we also consider the fully parametric case, where $g(s; \theta_3)$ is a prescribed parametric form, as described after (1.14), on the basis of which define $\hat{\theta}_3 = \arg\min_{\Theta_3} Q_{\log g}(\tilde{\theta}; \theta_3)$, and, with $\psi(s; \theta_3) = -(\partial/\partial s)g(s; \theta_3)/g(s; \theta_3)$,

$$\mathcal{J}_n(\theta, \sigma, \theta_3) = n^{-1}\sum_{t=1}^{n}\psi(E_t(\theta)/\sigma; \theta_3)^2,$$

$$r_i(\theta, \sigma, \theta_3) = \sum_{t=1}^{n}\psi(E_t(\theta)/\sigma; \theta_3)E'_{ti}(\theta), \qquad i = 1, 2,$$

and redefine $\hat{\theta}_i$, $i = 1, 2$, of (2.2) as

$$\hat{\theta}_i = \tilde{\theta}_i + \{R_i(\tilde{\theta})\mathcal{J}_n(\tilde{\theta}, \tilde{\sigma}, \hat{\theta}_3)\}^{-1}r_i(\tilde{\theta}, \tilde{\sigma}, \hat{\theta}_3), \qquad i = 1, 2.$$

We introduce the following additional assumptions.

ASSUMPTION A10. $\Theta_3$ is compact and $\theta_{03}$ is an interior point of $\Theta_3$.

ASSUMPTION A11. For all $\theta_3 \in \Theta - \{\theta_{03}\}$, $g(s; \theta_3) \neq g(s; \theta_{03})$ on a set of positive measure.

ASSUMPTION A12. In a neighborhood $\mathcal{N}$ of $\theta_{03}$, $\log g(s; \theta_3)$ is thrice continuously differentiable in $\theta_3$ for all $s$ and

$$\int_{-\infty}^{\infty}\left\{\sup_{\mathcal{N}}|g^{(k)}(s; \theta_3)| + \sup_{\mathcal{N}}|g^{(k,\ell)}(s; \theta_3)| + \sup_{\mathcal{N}}|g^{(k,\ell,m)}(s; \theta_3)|\right\}ds < \infty,$$

where $g^{(k)}$, $g^{(k,\ell)}$, $g^{(k,\ell,m)}$ represent partial derivatives of $g$ with respect to the $k$th, the $k$th and $\ell$th, and the $k$th, $\ell$th and $m$th elements of $\theta_3$, respectively.

ASSUMPTION A13. $\Omega_3 = E\{(\partial/\partial\theta_3)\log g(\varepsilon_t; \theta_{03})(\partial/\partial\theta_3^T)\log g(\varepsilon_0; \theta_{03})\}$ is positive definite.

THEOREM 2. *Let Assumptions* A1, A2(a), A3–A7 *and* A10–A13 *hold. Then as* $n \to \infty$, $n^{1/2}(\hat{\theta}_1 - \theta_{01})$, $D_n^{1/2}(\hat{\theta}_2 - \theta_{02})$ *and* $n^{1/2}(\hat{\theta}_3 - \theta_{03})$ *converge in distribution to independent* $N(0, \mathcal{J}^{-1}\Omega_1^{-1})$, $N(0, \mathcal{J}^{-1}\Omega_2^{-1})$ *and* $N(0, \Omega_3^{-1})$ *vectors, respectively, where the limiting covariance matrices are consistently estimated by* $\{\mathcal{J}_n(\tilde{\theta}, \tilde{\sigma}, \hat{\theta}_3)R_1(\tilde{\theta})/n\}^{-1}$, $\{\mathcal{J}_n(\tilde{\theta}, \tilde{\sigma}, \hat{\theta}_3)D_n^{-1}R_2(\tilde{\theta})D_n^{-1}\}^{-1}$ *and*

$$\left\{ n^{-1}\sum_{t=1}^{n}[(\partial/\partial\theta_3)\log g(E_t(\tilde{\theta})/\tilde{\sigma}; \hat{\theta}_3)][(\partial/\partial\theta_3^T)\log g(E_t(\tilde{\theta})/\tilde{\sigma}; \hat{\theta}_3)] \right\}^{-1},$$

*respectively.*

The proof (which entails an initial consistency proof for the implicitly defined extremum estimate $\hat{\theta}_3$) is omitted because it combines relatively standard arguments with elements of the proof of Theorem 1, notably concerning the $e_t(\theta_0) - \sigma_0\varepsilon_t$ issue. Our treatment of this would also lead to a theorem for $M$-estimates of $\theta_0$ minimizing (1.14) in which $\rho(s)$ is a completely specified function, not necessarily $\log g(s)$, but we omit this to conserve on space, and because the efficiency improvement of the paper's title would in general not be achieved.

Theorems 1 and 2 suggest locally more powerful (Wald-type) tests on $\theta_{01}$ than those implied by CLTs for Whittle estimates. For example, the hypothesis of short memory, $\xi_0 = 0$, can be efficiently tested, as can, say, the significance of AR coefficients in a FARIMA$(p_{11}, \xi_0, 0)$, for any unknown $\xi_0 > -\frac{1}{2}$. We can also efficiently investigate the question of relative success of deterministic and stochastic components in describing trending time series. For example, we can apply the theorems to test $\theta_{02} = 0$, or, with $p_2 = 1$, $p_2 = t^\tau$, test $\xi_0 = \tau + \frac{1}{2}$ against the one-sided alternative $\xi_0 > \tau + \frac{1}{2}$ [see the discussion after (1.9)]; in the first case rejection implies a significant deterministic trend, and in the latter, a dominant stochastic one. Tests based on $\hat{\theta}_2$ are in general more powerful than those based on least squares (see [31]) or generalized least squares (see [7]).

**4. Finite sample performance.** A small Monte Carlo study was carried out to investigate the success of our semiparametric estimates in small and moderate samples. Along with the value of $n$, major influential features seem likely to be the form of $g(s)$, the value of $\xi_0$ and the choice of $\phi$ and $L$.

We focused on the simple FARIMA$(0, \xi_0, 0)$ model for $y_t$ (knowing $\mu^T z_0 \equiv 0$) for:

(i) $\xi_0 = -0.25$ (antipersistent),

   (ii)  $\xi_0 = 0.25$ (stationary with long memory),
   (iii) $\xi_0 = 0.75$ (nonstationary but mean-reverting),
   (iv)  $\xi_0 = 1.25$ (nonstationary, non-mean-reverting).

For $\varepsilon_t$ we considered the following distributions [the scalings referred to producing var$(\varepsilon_t) = 1$]:

   (a)  $N(0, 1)$,
   (b)  $0.5N(-3, 1) + 0.5N(3, 1)$,
   (c)  (scaled) $0.05N(0, 25) + 0.95N(0, 1)$,
   (d)  (scaled) Laplace,
   (e)  (scaled) $t_5$.

These were mostly chosen for the sake of consistency with other Monte Carlo studies of adaptive estimates. The benchmark case (a), and the two (symmetric and asymmetric) mixed normal distributions (b) and (c), were used by Kreiss [19] in a stationary AR model, with kernel estimates of $\psi$, and by Newey [25] (in a cross-sectional regression model). Ling [22] used (b) in a FARIMA$(0, \xi_0, 0)$ model with kernel estimates of $\psi$. Kreiss [19] also used (d). The point of (e) is that it only just satisfies the minimal fourth moment condition on $\varepsilon_t$, Assumption A2(a). Kernel approaches, from [3] and [28] for location and regression models for independent observations, through Kreiss [19], Drost, Klaassen and Werker [8] and Koul and Schick [18] for short-memory time series models, and Hallin, Taniguchi, Serroukh and Choy [15], Hallin and Serroukh [14] and Ling [22] for long-memory ones, have been popular in the adaptive estimation literature. Besides requiring choice of a kernel and bandwidth (analogous to our $\phi$ and $L$), they typically involve one or more forms of trimming, in part due to the presence of a kernel density estimate in the denominator of the estimate of $\psi(s)$, and sometimes sample splitting and discretization of the initial estimate. Theorem 1 of course implies semiparametric efficient estimates using series estimation for short-memory models. For $\phi$ we used both (2.4) and (2.5), and tried $L = 1, 2, 3, 4$, with $n = 64$ and 128. For $\tilde{\xi} = \tilde{\theta}$ and $\tilde{\sigma}^2$ Velasco and Robinson's [29] estimates were employed, with a cosine bell taper; this is sufficient to satisfy Assumption A7 for all $\xi_0$ considered, albeit unnecessary when $\xi_0 = \pm 0.25$.

We report the Monte Carlo relative efficiency measure MSE$(\hat{\xi})$/MSE$(\tilde{\xi})$ (where $\hat{\xi} = \hat{\theta}$) for $L = 1$ on the basis of 1000 replications. Tables 1–5 present results for distributions (a)–(e), respectively, in case $n = 64$ only; generally asymptotic behavior was better approximated when $n = 128$. For $\varepsilon_t \sim N(0, 1)$, $\hat{\xi}$ is efficient when $\phi(s) = s$ for all $L \geq 1$, the efficiency improvement achieved in Table 1 for $L = 1$ being due to the tapering in $\hat{\xi}$; as anticipated, the unnecessarily complicated $\hat{\xi}$ based on larger $L$ makes matters somewhat worse. One expects relative efficiency to be roughly constant across $\xi_0$. The deviating results for $\xi_0 = -0.25$ and $\xi_0 = 1.25$ sometimes found in the tables are largely due to the following computational policy. The grid search to locate $\tilde{\xi}$ was confined to the interval $[-0.4, 1.75]$, and for the extreme $\xi_0$ some $\tilde{\xi}$ fell on the boundary (especially the lower one), while we correspondingly trimmed $\hat{\xi} < -0.4$ and $\hat{\xi} > 1.75$ to $\hat{\xi} = -0.4$ and $\hat{\xi} = 1.75$, respectively. This led to some underestimation of bias and variance, and consequent distortion of relative efficiency. However, there is considerable stability across $\xi_0$ in the symmetric mixed normal case (Table 2), and also small improvement with increasing $L$, but slight deterioration when $L = 4$ for the unbounded $\phi(s) = s$.

TABLE 1
$\varepsilon_t \sim N(0, 1)$

|  |  | $\phi(s) = s$ | | | | $\phi(s) = s(1 + s^2)^{-1/2}$ | | | |
|---|---|---|---|---|---|---|---|---|---|
|  | $L$ | 1 | 2 | 3 | 4 | 1 | 2 | 3 | 4 |
| $\xi_0$ | $-0.25$ | 0.62 | 0.62 | 0.62 | 0.62 | 0.66 | 0.67 | 0.63 | 0.65 |
|  | 0.25 | 0.47 | 0.48 | 0.51 | 0.61 | 0.49 | 0.52 | 0.53 | 0.60 |
|  | 0.75 | 0.46 | 0.49 | 0.53 | 0.62 | 0.50 | 0.54 | 0.55 | 0.60 |
|  | 1.25 | 0.47 | 0.50 | 0.52 | 0.61 | 0.52 | 0.53 | 0.52 | 0.56 |

For all tables, Monte Carlo MSE$(\hat{\xi})$/MSE$(\tilde{\xi})$ with $n = 64$ and 1000 replications.

TABLE 2
$\varepsilon_t \sim 0.5N(-3, 1) + 0.5N(3, 1)$

|  |  | $\phi(s) = s$ | | | | $\phi(s) = s(1 + s^2)^{-1/2}$ | | | |
|---|---|---|---|---|---|---|---|---|---|
|  | $L$ | 1 | 2 | 3 | 4 | 1 | 2 | 3 | 4 |
| $\xi_0$ | $-0.25$ | 0.92 | 0.92 | 0.83 | 0.90 | 0.94 | 0.93 | 0.82 | 0.83 |
|  | 0.25 | 0.90 | 0.91 | 0.89 | 0.93 | 0.91 | 0.91 | 0.88 | 0.89 |
|  | 0.75 | 0.90 | 0.91 | 0.89 | 0.94 | 0.90 | 0.92 | 0.89 | 0.89 |
|  | 1.25 | 0.88 | 0.89 | 0.88 | 0.92 | 0.89 | 0.89 | 0.87 | 0.87 |

TABLE 3
$\varepsilon_t \sim$ (scaled) $0.5N(0, 25) + 0.95N(0, 1)$

|  |  | $\phi(s) = s$ | | | | $\phi(s) = s(1 + s^2)^{-1/2}$ | | | |
|---|---|---|---|---|---|---|---|---|---|
|  | $L$ | 1 | 2 | 3 | 4 | 1 | 2 | 3 | 4 |
| $\xi_0$ | $-0.25$ | 0.71 | 0.71 | 0.62 | 0.77 | 0.81 | 0.76 | 0.63 | 0.70 |
|  | 0.25 | 0.84 | 0.76 | 0.65 | 0.74 | 0.77 | 0.67 | 0.60 | 0.54 |
|  | 0.75 | 0.85 | 0.79 | 0.70 | 0.79 | 0.80 | 0.78 | 0.69 | 0.63 |
|  | 1.25 | 1.01 | 0.96 | 0.81 | 0.82 | 0.91 | 0.83 | 0.74 | 0.68 |

TABLE 4
$\varepsilon_t \sim$ (scaled) Laplace

| | | $\phi(s) = s$ | | | | $\phi(s) = s(1+s^2)^{-1/2}$ | | | |
|---|---|---|---|---|---|---|---|---|---|
| | $L$ | 1 | 2 | 3 | 4 | 1 | 2 | 3 | 4 |
| $\xi_0$ | $-0.25$ | 1.07 | 0.85 | 0.92 | 0.96 | 1.04 | 0.90 | 0.60 | 0.61 |
| | 0.25 | 0.89 | 0.60 | 0.58 | 0.87 | 0.78 | 0.62 | 0.65 | 0.67 |
| | 0.75 | 0.56 | 0.52 | 0.55 | 0.81 | 0.51 | 0.53 | 0.53 | 0.54 |
| | 1.25 | 0.28 | 0.23 | 0.23 | 0.86 | 0.32 | 0.26 | 0.28 | 0.38 |

TABLE 5
$\varepsilon_t \sim$ (scaled) $t_5$

| | | $\phi(s) = s$ | | | | $\phi(s) = s(1+s^2)^{-1/2}$ | | | |
|---|---|---|---|---|---|---|---|---|---|
| | $L$ | 1 | 2 | 3 | 4 | 1 | 2 | 3 | 4 |
| $\xi_0$ | $-0.25$ | 0.58 | 0.54 | 0.53 | 0.65 | 0.55 | 0.53 | 0.55 | 0.60 |
| | 0.25 | 0.56 | 0.56 | 0.57 | 0.74 | 0.51 | 0.54 | 0.55 | 0.58 |
| | 0.75 | 0.58 | 0.58 | 0.62 | 0.75 | 0.51 | 0.56 | 0.57 | 0.61 |
| | 1.25 | 0.63 | 0.61 | 0.60 | 0.69 | 0.54 | 0.55 | 0.52 | 0.53 |

We find this also in the asymmetric mixed normal case (Table 3), though for the bounded $\phi(s) = s(1+s^2)^{-1/2}$, mainly the improvement continues to $L = 4$, and its magnitude, at each increase of $L$, is notable. For the Laplace distribution (Table 4) there is notable sensitivity to $\xi_0$, though increasing $L$ tends to improve efficiency, at least up to $L = 3$. For the $t_5$ distribution (Table 5) only small improvements, if any, were recorded after $L = 1$, as is not surprising for this small sample size, as asymptotic relative efficiency is 0.8; the deterioration with $\phi(s) = s$ at $L = 4$ is also not surprising due to the heavy tails. The results taken as a whole seem fairly encouraging, especially as the truncation (1.13) in computing residuals, which looms large in the theoretical component of this paper, would be expected to have some finite sample effect on $\hat{\xi}$ in our fractional setting.

**5. Final comments.** In various stationary, short-memory time series models, Kreiss [19], Drost, Klaassen and Werker [8], Koul and Schick [18] and others developed local asymptotic normality (LAN) and local asymptotic minimaxity (LAM) theory of Le Cam [20] and Hájek [12] to establish $\sqrt{n}$-consistent, asymptotically normal and asymptotically efficient estimates, and, further, adaptive estimates that achieve the same properties in the presence of nonparametric $g$. A similar approach was followed by Hallin et al. [15], Hallin and Serroukh [14] and Ling [22] in the case of stationary and nonstationary fractional models. LAN theory commences from a log-likelihood ratio, but in view of the difficulty in constructing likelihoods in a general non-Gaussian setting, the latter authors commenced not from the likelihood for $y_1, \ldots, y_n$ but from a "likelihood" for $y_1, \ldots, y_n$ and the infinite set of unobservable variables $\varepsilon_t$, $t \leq 0$, in terms of the density $g$ of $\varepsilon_t$, or a "conditional likelihood" for $y_1, \ldots, y_n$ given the $\varepsilon_t$, $t \leq 0$, or the $y_t$, $t \leq 0$. We do not employ such constructions and do not establish local optimality properties. However, the $M$-estimate efficiency bound we achieve is of course the same as the asymptotic variance resulting from a LAM/LAN approach.

Another motivation for our more elementary efficiency criterion is to allow space to focus on the main technical difficulty distinguishing asymptotic distribution theory for fractional models from that for short-memory ones. This is due to the need to approximate the truncated AR transforms $e_t = e_t(\theta_0)$ [see (1.13)] by scaled innovations $\sigma_0 \varepsilon_t$. Consider a simplified version of the problem in which $y_t = x_t$ a priori, so $\theta = \theta_1$, and define $\delta_t = e_t - \sigma_0 \varepsilon_t$. In the following section (relying heavily on Lemmas 13 and 14 of Section 7) we find that $E|\delta_t|^r \leq C t^{-r/2}$, $r \geq 2$, given a sufficient moment condition on $\varepsilon_t$. This property is useful in our proof that $e_t$ can be replaced by $\sigma_0 \varepsilon_t$ in $\hat{a}^{(L)}(E(\theta_0)/\sigma_0)$ (see Lemma 19). In some cases it is possible to show that the upper bound provides a sharp rate. Consider the stationary FARIMA$(0, \xi_0, 0)$ (cf. [14]), where $0 < \xi_0 = \zeta_0 < \frac{1}{2}$ and $x_t = v_t$, $t \in \mathbb{Z}$. Noting that $\operatorname{cov}(x_0, x_j) \geq j^{2\xi_0 - 1}/C$, $\alpha_j(\xi_0) \geq j^{-\xi_0 - 1}/C$ for $j > 0$,

$$
\begin{aligned}
E(\delta_t^2) &= \sum_{j=t}^{\infty} \sum_{k=t}^{\infty} \alpha_j(\xi_0) \alpha_k(\xi_0) \operatorname{cov}(x_j, x_k) \\
&\geq C^{-1} \sum_{\substack{j=t \\ 1 \leq |j-k| \leq t}}^{\infty} \sum_{k=t}^{\infty} j^{-\xi_0 - 1} k^{-\xi_0 - 1} |j-k|^{2\xi_0 - 1} \\
&\geq C^{-1} t^{2\xi_0 - 1} \sum_{j=t}^{\infty} \sum_{k=t+1}^{t+j} (jk)^{-\xi_0 - 1} \\
&\geq C^{-1} t^{2\xi_0 - 1} \sum_{j=t}^{2t} j^{-\xi_0} (t+j)^{-\xi_0 - 1} \\
&\geq (Ct)^{-1}.
\end{aligned}
$$

(This contrasts with the exponential rate occurring with ARMA models.) In this stationary FARIMA$(0, \xi_0, 0)$,

$$
\delta_t = \sum_{j=0}^{t-1} \alpha_j(\xi_0) x_{t-j} - \sigma_0 \varepsilon_t
$$

(5.1)
$$= \sum_{j=0}^{t-1} \alpha_j(\xi_0) v_{t-j} - \sigma_0 \varepsilon_t$$

$$= -\sum_{j=t}^{\infty} \alpha_{t+j}(\xi_0) v_{t-j}.$$

In our "asymptotically stationary" version of the FARIMA$(0, \xi_0, 0)$, also with $0 < \xi_0 < \frac{1}{2}$, we have $x_t = x_t^\#$, but again (5.1) results, from (1.4), (1.10), (1.11) and Lemma 5 of Section 7. In this connection, note that for general $\xi_0$, Ling [22] took $x_t = \Delta^{-m_0} v_t^\# + v_t \mathbb{1}(t \le 0)$ in place of our (1.2), but this different prescription of $x_t$ for $t \le 0$ makes no difference to $e_t$, which depends on $x_s$ for $s \ge 1$ only.

The above upper bound for $E|\delta_t|^r$, combined with the Schwarz inequality, is insufficient to deal completely with the replacement of $e_t$ by $\sigma_0 \varepsilon_t$, even when $\psi$ is smooth. Staying with the case $y_t = x_t$ a priori, the proofs of Theorems 1 and 2 entail establishing asymptotic normality of a quantity of the form $c_{1n} = n^{-1/2} \sum_{t=1}^n \psi(e_t) h_t$, where $h_t$ is $\{\varepsilon_s, s \le t-1\}$-measurable and has finite variance; $c_{1n}$ is called a "central sequence" by Hallin et al. [15] [see their (2.15) and (3.11)] and Hallin and Serroukh [14] [see their (2.4)]. Asymptotic normality of $c_{2n} = n^{-1/2} \sum_{t=1}^n \psi(\varepsilon_t) h_t$ follows straightforwardly from a martingale CLT. This leaves the relatively difficult task of showing that $c_{1n} - c_{2n} = o_p(1)$. In fact, our proof does not directly consider $c_{1n} - c_{2n}$ because we do not assume $\psi$ is smooth; we instead approximate the $e_t$ by the $\sigma_0 \varepsilon_t$ within the smooth estimate of $\psi$ and then appeal to mean square approximation of $\psi(\varepsilon_t)$ by its least-squares projection on the $\phi(\varepsilon_t)^\ell$, $\ell = 1, \ldots, L$, as $L \to \infty$, as in [25]. However, for this, $S_n = n^{-1/2} \sum_{t=1}^n \delta_t h_t$ [i.e., $c_{1n} - c_{2n}$ with $\psi(x)$ replaced by $x$] is relevant, and the sharper the bound we obtain for it the weaker some other conditions can be; we obtain $S_n = O_p((\log n)^{3/2} n^{-1/2})$.

The same kind of issue arises in theory for Whittle estimation. For short-memory stationary processes, with $\xi_0 = 0$, Hannan [16] established the CLT for various Whittle estimates. His proof does not work under stationary long memory, $0 < \xi_0 < \frac{1}{2}$, due to the bad behavior of the periodogram and spectral density at low frequencies. However, in this case Fox and Taqqu [9], Dahlhaus [6] and Giraitis and Surgailis [11] delicately exploited a kind of balance between these quantities in order to establish CLTs. The CSS estimate minimizing $\sum_{t=1}^n e_t^2(\theta)$ [see Remark 5 in Section 3 concerning (1.14)] is not one of those considered by these authors, but its CLT requires showing $S_n = o_p(1)$, which entails similar challenge to results they established for the somewhat different quadratic forms arising from their parameter estimates. Our results for replacing $e_t$ by $\sigma_0 \varepsilon_t$ can be employed to provide a proof of asymptotic normality of the CSS version of the Whittle estimate. Whittle and adaptive estimation are both areas in which asymptotic results are qualitatively the same across short and long memory, but sufficient methods of proof significantly differ.

**6. Proof of Theorem 1.** The consistency of the covariance matrix estimates is implied by the proof of the CLT. By far the most significant features of this are accomplished in the lemmas in the following section. Their application is mostly relatively straightforward, and is thus described here in abbreviated form. For notational convenience we now write $\theta_3 = \sigma$ and augment $\theta$ as $\theta = (\theta_1^T, \theta_2^T, \theta_3)^T$. We also abbreviate $\sum_{t=1}^{n}$ to $\sum_t$, and $E_t(\theta_0)$, $E(\theta_0)$, $E_{ti}(\theta_0)$ to $E_t$, $E$, $E_{ti}$, respectively, $i = 1, 2$. By the mean value theorem, for $i = 1, 2$,

$$\hat{\theta}_i - \theta_{0i} = \left\{ I_{p_i} + \frac{R_i(\tilde{\theta})^{-1}}{\mathcal{J}_L(\bar{\theta})} \bar{S}_{Lii} \right\} (\tilde{\theta}_i - \theta_{0i})$$
$$+ \frac{R_i(\tilde{\theta})^{-1}}{\mathcal{J}_L(\bar{\theta})} \left\{ \sum_{j=1, j \neq i}^{3} \bar{S}_{Lij}(\tilde{\theta}_j - \theta_{0j}) + r_{Li}(\theta_0) \right\},$$

where, with $[S_{Li1}(\theta), S_{Li2}(\theta), S_{Li3}(\theta)] = (\partial/\partial\theta^T) r_{Li}(\theta)$, each row of $\bar{S}_{Lij}$ is formed from the corresponding row of $S_{Lij}(\theta)$ by replacing $\theta$ by $\bar{\theta}$ such that $\|\bar{\theta} - \theta_0\| \leq \|\tilde{\theta} - \theta_0\|$ where $\|A\| = \{\mathrm{tr}(A^T A)\}^{1/2}$. Write $D_{1n} = D_{3n} = n^{1/2}$, $D_{2n} = D_n$ and define $\mathcal{N} = \{\theta : \|D_{in}(\theta_i - \theta_{0i})\| \leq 1, \ i = 1, 2, 3\}$. The result follows if

$$(6.1) \qquad \sup_{\mathcal{N}} \|D_{in}^{-1} \{R_i(\theta) - R_i(\theta_0)\} D_{in}^{-1}\| \xrightarrow{p} 0, \qquad i = 1, 2,$$

$$(6.2) \quad \sup_{\mathcal{N}} \|D_{in}^{-1} \{S_{Lij}(\theta) - S_{Lij}(\theta_0)\} D_{jn}^{-1}\| \xrightarrow{p} 0, \qquad i = 1, 2, j = 1, 2, 3,$$

$$(6.3) \qquad \sup_{\mathcal{N}} |\mathcal{J}_L(\theta) - \mathcal{J}_L(\theta_0)| \xrightarrow{p} 0,$$

$$(6.4) \qquad D_{in}^{-1} R_i(\theta_0) D_{in}^{-1} \xrightarrow{p} \Omega_i, \qquad i = 1, 2,$$

$$(6.5) \qquad \{R_i(\theta_0)\mathcal{J}_L(\theta_0)\}^{-1} S_{Lij}(\theta_0) \xrightarrow{p} -I_{p_i}\mathbb{1}(i = j),$$
$$\qquad\qquad\qquad\qquad\qquad\qquad\qquad\qquad i = 1, 2, j = 1, 2, 3,$$

$$(6.6) \qquad \mathcal{J}_L(\theta_0) \xrightarrow{p} \mathcal{J},$$

$$(6.7) \qquad \begin{bmatrix} n^{-1/2} r_1 \\ D_n^{-1} r_2 \end{bmatrix} \xrightarrow{d} N\left(0, \begin{bmatrix} \mathcal{J}\Omega_1 & 0 \\ 0 & \mathcal{J}\Omega_2 \end{bmatrix}\right),$$

$$(6.8) \qquad D_{in}^{-1} \{r_{Li}(\theta_0) - r_i\} \xrightarrow{p} 0, \qquad i = 1, 2,$$

where

$$r_1 = \sum_t \psi(\varepsilon_t) \varepsilon'_{t1}, \qquad r_2 = \sum_t \psi(\varepsilon_t) E'_{t2},$$

with $\varepsilon'_{t1} = (\partial/\partial\theta_1^{(+)T})\alpha(B;\theta_1^{(+)})/\sigma_0 = \gamma(B;\nu_0)\varepsilon_t$.

The most difficult and distinctive problems occur in (6.8) for $i = 1$, which faces the $e_t - \sigma_0\varepsilon_t$ problem, as well as the increasing $L$, in the presence of normalization only by $D_{1n}^{-1}$. The first of these aspects is also in (6.1) and (6.4), and both are in (6.2), (6.3), (6.5) and (6.6), but the normalizations make (6.4)–(6.6) much easier to deal with and the proof details are otherwise relatively standard, albeit lengthy. The same may also be said for (6.1)–(6.3), except for the approximation of the fractional difference $\Delta^{\xi_0}$ by $\Delta^\xi$ for $|\xi - \xi_0| \le n^{-1/2}$, bearing in mind that "nonstationary" values of $\xi$, $\xi_0$ are permitted. The basic steps in proving (6.1)–(6.3) are illustrated by the least complicated case (6.1). By elementary inequalities it suffices to show that $\sup_{\mathcal{N}}\sum_t\|D_{in}^{-1}(e'_{ti}(\theta) - e'_{ti}(\theta_0))\|^2 \overset{p}{\to} 0$, $i = 1, 2$. Write $\alpha = \alpha(B;\theta^{(-)})$, $\alpha' = \alpha'(B;\theta^{(-)})$ with $\alpha_0, \alpha'_0$ denoting these quantities at $\nu = \nu_0$. For $i = 2$, it suffices to apply Lemmas 1, 2, 3 and (with $m = \xi_0$) 4, the $j$th elements of $\alpha_0(\Delta^\xi - \Delta^{\xi_0})z_{2t}$ and $(\alpha - \alpha_0)\Delta^{\xi_0}z_{2t}$ being, respectively, $O(n^{-1/2}(\log t)t^{\chi_j - \xi_0})$ and $O(n^{-1/2}t^{\chi_j - \xi_0})$ uniformly in $\mathcal{N}$, noting that $\xi_0 > -\frac{1}{2}$ and $\chi_j \ge \xi_0 - \frac{1}{2}$ implies $\chi_j > -1$ and $\xi_0 < \chi_j + 1$. For $i = 1$, the terms in $z_{2t}$ are dealt with similarly, while Lemmas 1–4 give, for example, $\alpha'_0(\Delta^\xi - \Delta^{\xi_0})(s_t + \mu^*t^{\xi_0}) = O(n^{-1/2}(\log t)^2)$ and $(\alpha' - \alpha'_0)\Delta^{\xi_0}(s_t + \mu^*t^{\xi_0}) = O(n^{-1/2})$ uniformly in $\mathcal{N}$. In the above we apply first Lemma 3, then Lemma 1 and then Lemma 2, noting that in case (ii) of Lemma 1 must be used (either for a leading term or remainder) the coefficient of $s^j$ in the expansion of $-\log(1-s)$, and thus of $(-\log(1-s))^r$, is positive for all $j \ge 1$, so for nonnegative sequences $g_t, h_t$, such that $g_t \le h_t$, we have $|(-\log\Delta)^r g_t| \le |(-\log\Delta)^r x_t|$. So far as contributions from $x_t$ are concerned, from Lemma 5

$$\sup_{\mathcal{N}}\|(\alpha' - \alpha'_0)\Delta^{\xi_0}x_t\| \le \sum_{j=0}^{t-1}\left\{\sup_{\mathcal{N}}\|\alpha'_j - \alpha'_{0j}\|\right\}\{|\Delta^{\xi_0}v^\#_{t-j}| + |(\log\Delta)\Delta^{\xi_0}v^\#_{t-j}|\},$$

where $\alpha'_j$, $\alpha'_{0j}$ are the $j$th Fourier coefficients of $\alpha'$, $\alpha'_0$. By the mean value theorem and Lemma 6 this has second moment $O(n^{-1})$. The same result holds for $\alpha'_0(\Delta^\xi - \Delta^{\xi_0})x_t$ after taking $m = m_0$ in Lemma 4, noting that its supremum over $\mathcal{N}$ is bounded by

$$Cn^{-1/2}\|\alpha'_0\Delta^{\xi_0}x_t\| + Cn^{-1/2}\|(\log\Delta)\alpha'_0\Delta^{\xi_0}x_t\| + Cn^{-1}\left(\sum_{j=1}^t v^2_{t-j}\right)^{1/2}$$

and applying Lemmas 5 and 6. The proof of (6.1) is readily completed.

Before coming to (6.8), we briefly discuss (6.7). Consider variates $U = (n^{-1/2}r_1^T, (D_{2n}^{-1}r_2)^T)^T$, $V = \lambda^T(EUU^T)^{-1/2}U$ for a $(p_1 + p_2) \times 1$ vector $\lambda$ such that $\lambda^T\lambda = 1$. We have $EV = 0$, $EV^2 = 1$, since $E\psi(\varepsilon_0) = 0$ and $\varepsilon'_{t1}$ is independent of $\varepsilon_t$, so (6.7) follows from Theorem 2 of [27] if

$$(6.9) \qquad \sum_t \begin{bmatrix} n^{-1/2}\varepsilon'_{t1} \\ D_n^{-1}E'_{t2} \end{bmatrix} \begin{bmatrix} n^{-1/2}\varepsilon'_{t1} \\ D_n^{-1}E'_{t2} \end{bmatrix}^T \xrightarrow{p} \begin{bmatrix} \Omega_1 & 0 \\ 0 & \Omega_2 \end{bmatrix},$$

$$\sum_t \psi(\varepsilon_t^2)\{n^{-1}\|\varepsilon'_{t1}\|^2 \mathbb{1}(\|\psi(\varepsilon_t)\varepsilon'_{t1}\| \geq \delta n^{1/2})$$

$$(6.10) \qquad\qquad + \|D_n^{-1}E'_{t2}\|^2 \mathbb{1}(\|\psi(\varepsilon_t)D_n^{-1}E'_{t2}\| \geq \delta)\} \xrightarrow{p} 0$$

for any $\delta > 0$. The proof of (6.9) follows from Lemmas 1 and 3 and approximating sums by integrals, while that of (6.10) follows from stationarity and finite variance of $\psi(\varepsilon_t)$ and $\varepsilon'_{t1}$ and the slowly changing character of $z_{2t}$.

We prove (6.8) only for $i = 1$, the case $i = 2$ involving some of the same steps but being much easier. Define $\Xi^{(L)}(s) = \phi^{(L)}(s) - E\phi^{(L)}(\varepsilon_t)$, $W^{(L)} = E\{\Xi^{(L)}(\varepsilon_t)\Xi^{(L)}(\varepsilon_t)^T\}$. It follows from Lemma 8 that $W^{(L)}$ is nonsingular, and thence we define $a^{(L)} = W^{(L)-1}w^{(L)}$ where $w^{(L)} = E\{\phi'^{(L)}(\varepsilon_t)\} = E\{\phi^{(L)}(\varepsilon_t)\psi(\varepsilon_t)\}$, by integration by parts, as in [2] and as justified under our conditions by Lemma 2.2 of [25]. Defining also $\bar{\psi}^{(L)}(\varepsilon_t; a^{(L)}) = a^{(L)T}\Xi^{(L)}(\varepsilon_t)$, we have

$$n^{-1/2}\{r_{L1}(\theta_0) - r_1\} = \sum_{i=1}^{4}\sum_{j=1}^{2} A_{ij} - A_{11},$$

where

$$A_{ij} = n^{-1/2}\sum_t B_{it}C_{jt}$$

and

$$B_{1t} = \psi(\varepsilon_t),$$
$$B_{2t} = \bar{\psi}^{(L)}(\varepsilon_t; a^{(L)}) - \psi(\varepsilon_t),$$
$$B_{3t} = \psi^{(L)}(\varepsilon_t; \hat{a}^{(L)}(\varepsilon)) - \bar{\psi}^{(L)}(\varepsilon_t; a^{(L)}),$$
$$B_{4t} = \tilde{\psi}_t^{(L)}(\theta_0, \sigma_0) - \psi^{(L)}(\varepsilon_t; \hat{a}^{(L)}(\varepsilon)),$$
$$C_{1t} = \sigma_0\varepsilon'_{t1}, \qquad C_{2t} = E'_{t1} - \sigma_0\varepsilon'_{t1}.$$

Since $\varepsilon'_{t1}$ is $\{\varepsilon_s, s < t\}$-measurable and $E\|\varepsilon'_{01}\|^2 \leq C\|\Omega_1\| < \infty$, while $B_{2t}$ has zero mean, $E\|A_{21}\|^2 \leq CEB_{20}^2 \to 0$ as $L \to \infty$ from [10], pages 74–77, and [25], Lemma 2.2, since the moments of $\phi(\varepsilon_0)$ characterize its distribution under Assumptions A2 and A8.

Before discussing other $A_{ij}$ define

$$\mu_a = 1 + E\{|\varepsilon_t|^a \mathbb{1}(|\varepsilon_t| > 1)\},$$

for $a > 0$, and the following sequences:

$$
\begin{aligned}
\rho_{aL} &= CL & &\text{if } a = 0, \\
&= (CL)^{aL/\omega} & &\text{if } a > 0 \text{ and Assumption A2(b) holds,} \\
&= C^L & &\text{if } a > 0 \text{ and Assumption A2(c) holds,}
\end{aligned}
$$

suppressing reference in $\rho_{aL}$ to the arbitrarily large constant $C$; and also

$$
\pi_L = (\log L)\eta^{2L}\mathbb{1}(\varphi < 1) + (L\log L)\eta^{2L}\mathbb{1}(\varphi = 1) + (\log L)(\eta\varphi)^{2L}\mathbb{1}(\varphi > 1),
$$

for $L > 1$.

Write $A_{31} = (b_{1n} - b_{2n}b_{3n})\{\hat{a}^{(L)}(\varepsilon) - a^{(L)}\} - b_{2n}b_{3n}a^{(L)}$, where $b_{1n} = n^{-1/2}\sigma_0\sum_t \varepsilon'_{t1}\Xi^{(L)}(\varepsilon_t)^T$, $b_{2n} = n^{-1}\sum_t \varepsilon'_{t1}$, $b_{3n} = n^{-1/2}\sigma_0\sum_t \Xi^{(L)}(\varepsilon_t)^T$. We have $E|\phi(\varepsilon_0)|^r \leq \mu_{\kappa r}$, and thus from Lemma 9

$$
E\|b_{1n}\|^2 + E\|b_{3n}\|^2 \leq C\sum_{\ell=1}^{L}(E\|\varepsilon'_{01}\|^2 + 1)E\phi^{2\ell}(\varepsilon_0) \leq \rho_{2\kappa L}.
$$

Since $b_{2n} = O_p(n^{-1/2}\log n)$ from Lemma 17, we deduce from Lemma 10 that

$$
(6.10) \qquad A_{31} = O\left(\frac{L\rho_{2\kappa L}\pi_L}{n^{1/2}}(\log n + L^{1/2}\rho_{4\kappa L}^{1/2}\pi_L)\right).
$$

Before imposing Assumption A9, we estimate $A_{41}$, which can be written

$$
(6.11) \qquad n^{-1/2}\sigma_0\left[\sum_t \varepsilon'_{t1}\{\Phi^{(L)}(E_t/\sigma_0) - \Phi^{(L)}(\varepsilon_t)\}\right]\hat{a}^{(L)}(E/\sigma_0)
$$

$$
(6.12) \qquad + n^{-1/2}\sigma_0\sum_t \varepsilon'_{t1}\Phi^{(L)}(\varepsilon_t)^T\{\hat{a}^{(L)}(E/\sigma_0) - \hat{a}^{(L)}(\varepsilon)\}.
$$

The square-bracketed quantity in (6.11) has norm bounded by

$$
(6.13) \qquad \left(\sum_{\ell=1}^{L}\left\|\sum_t \varepsilon'_{ti}\delta_{\ell t}\right\|^2\right)^{1/2} + n^{-1}\left\|\sum_t \varepsilon'_{ti}\right\|\left\{\sum_{\ell=1}^{L}\left(\sum_t \delta_{\ell t}\right)^2\right\}^{1/2},
$$

where $\delta_{\ell t} = \phi_\ell(E_t/\sigma_0) - \phi_\ell(\varepsilon_t)$. We have

$$
(6.14) \qquad \delta_{\ell t} = \phi'_\ell(\varepsilon_t)d_t + \tfrac{1}{2}\phi''_\ell(\bar{\varepsilon}_t)d_t^2,
$$

where $|\bar{\varepsilon}_t - \varepsilon_t| \leq |d_t|$, $d_t = E_t/\sigma_0 - \varepsilon_t$. Now $e_t = \alpha(B;\theta_{01})(s_t + \mu^* t^{\xi_0} + x_t)$, and from Lemma 5 [see also (1.13)]

$$
\alpha(B;\theta_{01})x_t = \alpha(B;\theta_{01}^{(+)})v_t^{\#} = \sigma_0\varepsilon_t - \sum_{j=0}^{\infty}\alpha_{t+j}(\theta_{01}^{(+)})v_{-j} = \sigma_0\varepsilon_t + d_{1t},
$$

where

$$
d_{1t} = -\sum_{j=1}^{\infty}\lambda_{jt}\varepsilon_{t-j}, \qquad \lambda_{jt} = \sum_{k=0}^{j}\alpha_{k+t}(\theta_{01}^{(+)})\beta_{j-k}(\theta_{01}^{(+)}),
$$

where $\beta_j(\theta_{01}^{(+)})$ is the coefficient of $s^j$ in $\alpha(s;\theta_{01}^{(+)})^{-1}$. Since $\alpha(B;\theta_0)s_t = o(t^{-1/2})$ and $\alpha(B;\theta_{01})t^{\xi_0} = \alpha(1;\theta_0^{(-)})\Gamma(\xi_0+1) + O(t^{-1})$ from Lemma 1, it follows that

$$d_t = d_{1t} + d_2 + d_3 + o(t^{-1/2}), \tag{6.15}$$

where $d_2 = n^{-1}\sum_{j=0}^{\infty}(\sum_t \lambda_{jt})\varepsilon_{-j}$, $d_3 = n^{-1}\sum_t \varepsilon_t$. From Lemmas 13, 14 and 18, for $2 \le r \le 4$ under Assumption A2(a) and $r > 4$ under Assumptions A2(b) and A2(c),

$$E|d_{1t}|^r \le (Cr)^{2r}t^{-r/2}\mu_{r_+}^{r/r_+}, \tag{6.16}$$

$$E|d_2|^r + E|d_3|^r \le (Cr)^{2r}n^{-r/2}\mu_{r_+}^{r/r_+}, \tag{6.17}$$

where $r_+$ is the smallest even integer such that $r \le r_+$. Returning to (6.13), we have

$$\left\|\sum_t \varepsilon_{t1}'\delta_{\ell t}\right\| \le \left\|\sum_t \varepsilon_{t1}'\{\phi_\ell'(\varepsilon_t) - E\phi_\ell'(\varepsilon_0)\}d_{1t}\right\| \tag{6.18}$$

$$+ \left\|\sum_t \varepsilon_{t1}'\{\phi_\ell'(\varepsilon_t) - E\phi_\ell'(\varepsilon_0)\}\right\|(|d_2| + |d_3|) \tag{6.19}$$

$$+ |E\phi_\ell'(\varepsilon_0)|\left\|\sum_t \varepsilon_{t1}'d_{1t}\right\| \tag{6.20}$$

$$+ |E\phi_\ell'(\varepsilon_0)|\left\|\sum_t \varepsilon_{t1}'\right\|(|d_2| + |d_3|) \tag{6.21}$$

$$+ \left\|\sum_t \varepsilon_{t1}'\phi_\ell''(\bar{\varepsilon}_t)d_t^2\right\|. \tag{6.22}$$

Now

$$\begin{aligned}|\phi_\ell'(s)| &= \ell|\phi'(s)\phi^{\ell-1}(s)| \\ &\le C\ell(1 + |\phi(s)|^K)\{\mathbb{1}(|s| \le 1) + |s|^{\kappa(\ell-1)}\mathbb{1}(|s| > 1)\} \\ &\le C\ell\{\mathbb{1}(|s| \le 1) + |s|^{\kappa(\ell-1+K)}\mathbb{1}(|s| > 1)\},\end{aligned} \tag{6.23}$$

and since $\varepsilon_t$ is independent of $\varepsilon_{ti}'d_{1t}$, the right-hand side of (6.18) is

$$O_p\left(\{E\phi_\ell'(\varepsilon_0)^2\}^{1/2}\sum_t (E\|\varepsilon_{ti}'\|^4 Ed_{1t}^4)^{1/2}\right) = O_p(\ell\mu_{2\kappa(\ell+K)}^{1/2}\log n),$$

using (6.16). The same bound applies to (6.19)–(6.21), proceeding similarly and using respectively (6.17), Lemma 16, and (6.17) with Lemma 17; note that it is the second factor in (6.20) which leads to the main work in handling the quantity $S_n$ discussed in Section 5. So far as (6.22) is concerned, note that as in (6.23),

$$|\phi_\ell''(s)| \le C\ell^2 \{\mathbb{1}(|s| \le 1) + |s|^{\kappa(\ell-1+2K)} \mathbb{1}(|s| > 1)\},$$

so by the $c_r$-inequality ([23], page 157) (6.22) is bounded by

$$(6.24) \quad C^{\kappa\ell+1}\ell^2 \sum_t \|\varepsilon_{t1}'\| \{d_{1t}^2 + d_{1t}^2 |\varepsilon_t|^{\kappa(\ell+K)} + |d_{1t}|^{\kappa(\ell+K)+2}\}$$

$$+ C^{\kappa\ell+1}\ell^2 \sum_t \|\varepsilon_{t1}'\| \{(d_t - d_{1t})^2 (1 + |\varepsilon_t|^{\kappa(\ell+K)})$$

$$(6.25) \quad + |d_t - d_{1t}|^{\kappa(\ell+K)+2}\}.$$

By (6.16) and Hölder's and Jensen's inequalities, (6.24) has expectation bounded by

$$C^{\kappa\ell+1}\ell^2 \left\{ \mu_{\kappa(\ell+K)} \log n + \sum_t (E|d_{1t}|^{2\kappa(\ell+K)+4})^{1/2} \right\}$$

$$\le C(C\ell)^{2\kappa\ell}\ell^2 \mu_{r_\ell}^{1/2} \log n,$$

$r_\ell$ being the smallest integer such that $r_\ell \ge 2\kappa(\ell+K) + 4$. From (6.14) and (6.17), (6.25) is of smaller order in probability. It follows from Lemma 9 that

$$\left( \sum_{\ell=1}^{L} \left\| \sum_t \varepsilon_{t1}' \delta_{\ell t} \right\|^2 \right)^{1/2} = O_p((CL)^{2\kappa L+2} \rho_{2\kappa L}^{1/2} \log n).$$

By a similar but easier proof, the second term in (6.13) has the same bound, and by Lemmas 10 and 19,

$$(6.11) = O_p((CL)^{2\kappa L+3} \rho_{2\kappa L} \pi_L n^{-1/2} \log n).$$

Next, from similar but simpler arguments to those above,

$$n^{-1/2} \left\| \sum_t \varepsilon_{t1}' \Phi^{(L)}(\varepsilon_t)^T \right\| = O_p(\rho_{2\kappa L}^{1/2} \log n).$$

Application of Lemma 9 indicates that (6.12) is

$$O_p(\rho_{2\kappa L}^2 \pi_L^2 (L^2 n^{-1/2} \log n + (CL)^{4\kappa L+3} n^{-1} (\log n)^2)).$$

Thus

$$(6.26) \quad A_{41} = O_p(\rho_{2\kappa L} \pi_L (\rho_{2\kappa L} \pi_L L^2 + (CL)^{2\kappa L+3}$$

$$+ \rho_{2\kappa L} \pi_L (CL)^{4\kappa L+3} n^{-1/2}) n^{-1/2} \log n).$$

Comparison of (6.10) and (6.26) indicates that $A_{31}$ is dominated by $A_{41}$, whose behavior under Assumption A9 we thus now consider. Take $\kappa = 0$. From Lemma 9, under Assumption A9(a)

$$A_{41} = O_p(L^4 \pi_L^2 n^{-1/2} \log n)$$
$$= O_p\left(\exp\left[\log n\left\{\frac{4 \log L + \log \log n + 2 \log \pi_L}{\log n} - \frac{1}{2}\right\}\right]\right),$$

which is $o_p(1)$ if $\limsup \log \pi_L / \log n < \frac{1}{4}$, as is clearly implied by (3.4). Now take $\kappa > 0$ under Assumption A2(b). From Lemma 9, under Assumption A9(b)

$$A_{41} = O_p((L^{4\kappa L/\omega + 2} + L^{2\kappa L(1 + 1/\omega) + 3})n^{-1/2} \log n) = o_p(1),$$

on proceeding as before. Under Assumption A2(c), Lemma 9 and Assumption A9(c) give

$$A_{41} = O_p((CL)^{2\kappa L} n^{-1/2} \log n) = o_p(1).$$

To consider $A_{12}$, we can proceed as earlier to write

$$E'_{t1} - \varepsilon'_{t1} = D_{1t} + D_2 + D_3 + (t^{-1/2} \log t),$$

where

$$D_{1t} = -\sum_{j=1}^{\infty} \tilde{\lambda}_t \varepsilon_{t-j}, \qquad D_2 = n^{-1} \sum_{j=0}^{\infty} \left(\sum_t \tilde{\lambda}_{jt}\right) \varepsilon_{-j}, \qquad D_3 = n^{-1} \sum_t \varepsilon'_{t1},$$

and $\tilde{\lambda}_{jt} = \sum_{k=0}^{j} (\partial/\partial \theta_1^{(+)T}) \alpha_{k+t}(\theta_{01}^{(+)}) \beta_{j-k}(\theta_{01}^{(+)})$. Using (7.23) and (7.24) of Lemma 13, we deduce that $|\tilde{\lambda}_{jt}| \leq C(\log t) j^{\zeta_0} t^{-\zeta_0 - 1}$, $j \leq t$, and $|\tilde{\lambda}_{jt}| \leq C(\log t) j^{\zeta_0 - 1} \max(j^{-\zeta_0}, t^{-\zeta_0})$, $j > t$, and then proceeding as in Lemma 14, that $\sum_{j=0}^{\infty} \tilde{\lambda}_{jt}^2 \leq Ct^{-1} \log^2 t$, $\sum_{j=0}^{\infty} (\sum_{t=1}^{n} \tilde{\lambda}_{jt})^2 \leq Cn \log^2 n$. Noting that $E(\sum_t \psi(\varepsilon_t) \times D_{1t})^2 \leq C \sum_t ED_{1t}^2$, using also Lemma 17 and proceeding as in the proof for (6.11), it follows that $A_{12} = O_p(n^{-1/2} \log^{3/2} n)$.

The remainder of the proof of (6.8) with $i = 1$ deals in similar if easier ways with quantities already introduced and is thus omitted. □

## 7. Technical lemmas.

To simplify lemma statements, we take it for granted that, where needed, Assumptions A1–A9 hold.

Part (ii) of the following lemma is only needed to show that $s_t$ in (1.9) contributes negligibly, in particular when it includes $\tau_1 \leq \xi_0 - 1$.

LEMMA 1. (i) For $w_t = t^{\gamma}$ with $\gamma > -1$ and $\xi \in (-\frac{1}{2}, \gamma + 1)$,

$$\Delta^{\xi} w_t^{\#} = \frac{\Gamma(\gamma + 1)}{\Gamma(\gamma - \xi + 1)} t^{\gamma - \xi} + O(t^{\gamma - \xi - 1} + t^{\gamma - m - 1} \mathbb{1}(\xi > 0)),$$

*as $t \to \infty$, where $m$ is the integer such that $\xi - 1 < m \leq \xi$.*

(ii) *For $w_t = (\log t)^r t^\gamma$, $r \geq 0$, $\xi > -\frac{1}{2}$,*

$$\Delta^\xi w_t^\# = O(t^{\max(\gamma, -1) - \xi + d}) \qquad as \ t \to \infty,$$

*for any $\delta > 0$.*

PROOF. (i) The proof when $\xi$ is a nonnegative integer is straightforward, so we assume this is not the case. We have

$$(7.1) \qquad \sum_{j=0}^{\infty} j^k \Delta_j(-\xi) = 0, \qquad j = 0, \dots, m,$$

when $m \geq 0$ and $\xi > 0$, $(1-s)^\xi$ and its first $m$ derivatives in $s$ being zero at $s = 1$. With $a_k = \Delta_k(-\gamma)$,

$$
\begin{aligned}
\Delta^\xi w_t^\# &= \sum_{j=0}^{t-1} \Delta_j(-\xi)(t-j)^\gamma \\
&= t^\gamma \sum_{j=0}^{t-1} \Delta_j(-\xi) \sum_{k=0}^{\infty} a_k (j/t)^k \\
&= -t^\gamma \sum_{k}{}' (t-k)^{-k} a_k \sum_{j=t}^{\infty} j^k \Delta_j(-\xi) \mathbb{1}(m \geq 0) \\
&\quad + t^\gamma \sum_{k}{}'' (t-k)^{-k} a_k \sum_{j=0}^{t-1} j^k \Delta_j(-\xi),
\end{aligned}
$$

where $\sum_{k}{}' = \sum_{k=0}^{m}$, $\sum_{k}{}'' = \sum_{k=\max(m+1,0)}^{\infty}$ and we apply (7.1). By Stirling's approximation

$$(7.3) \qquad \left| \Delta_j(-\xi) - \frac{j^{-\xi-1}}{\Gamma(-\xi)} \right| \leq C j^{-\xi-2}, \qquad j \geq 1,$$

so (7.2) differs from

$$(7.4) \qquad \frac{t^\gamma}{\Gamma(-\xi)} \left\{ -\sum_{k}{}' (t-k)^{-k} a_k \sum_{j=t}^{\infty} j^{k-\xi-1} \mathbb{1}(m \geq 0) \right.$$
$$\left. + \sum_{k}{}'' (t-k)^{-k} a_k \sum_{j=0}^{t-1} j^{k-\xi-1} \right\}$$

by

$$(7.5) \qquad O\left(t^\gamma \sideset{}{'}\sum_k t^{-k}|a_k| \sum_{j=t}^\infty j^{k-\xi-2} \mathbb{1}(m \geq 0) + t^\gamma \sideset{}{''}\sum_k t^{-k}|a_k| \sum_{j=0}^{t-1} j^{k-\xi-2}\right).$$

Now

$$
\begin{aligned}
(7.6) \qquad
&\sum_{j=0}^{t-1} j^{-\alpha} = t^{1-\alpha}/(1-\alpha) + O(t^{-\alpha}), \qquad \alpha < 1, \\
&\sum_{j=t}^\infty j^{-\alpha} = t^{1-\alpha}/(\alpha-1) + O(t^{-\alpha}), \qquad \alpha > 1.
\end{aligned}
$$

Thus (7.5) is

$$
\begin{aligned}
O\Bigg(&t^{\gamma-\xi-1} \sideset{}{'}\sum_k \frac{|a_k|}{\xi+1-k} \mathbb{1}(m \geq 0) \\
&+ t^{\gamma-m-1}\frac{|a_{m+1}|}{\xi-m}\mathbb{1}(m \geq -1) + t^{\gamma-\xi-1}\sideset{}{'''}\sum_k \frac{|a_k|}{k-\xi-1}\Bigg) \\
={}& O\Bigg(t^{\gamma-\xi-1}\bigg\{\sideset{}{'}\sum_k |a_k|\mathbb{1}(m \geq 0) \\
&\qquad + \sideset{}{''}\sum_k (k+1)^{-1}|a_k| + Ct^{\gamma-m-1}|a_{m+1}|\mathbb{1}(m \geq -1)\bigg\}\Bigg),
\end{aligned}
$$

where $\sum_k''' = \sum_{k=\max(m+2,0)}^\infty$. The first sum in braces is finite because $m$ and the $a_k$ are, while the second sum is finite because $|a_k| \leq Ck^{-\gamma-1}$. Thus since $\gamma > -1$, (7.5) is $O(t^{\gamma-m-1})$ for $\xi > 0$ and $O(t^{\gamma-\xi})$ for $\xi < 0$. Applying (7.6) again, (7.4) is

$$\frac{t^{\gamma-\xi}}{\Gamma(-\xi)} \sum_{k=0}^\infty \frac{a_k}{k-\xi} + O(t^{\gamma-\xi-1}),$$

and the leading term is $\{\Gamma(\gamma+1)/\Gamma(\gamma-\xi+1)\}t^{\gamma-\xi}$, from [30], page 260.

(ii) We have

$$\Delta^\xi w_t^\# = \sum_{j=0}^{t-1} \Delta_j(-\xi)\{\log(t-j)\}^r (t-j)^\gamma.$$

Noting that $\Delta_j(-\xi) = O(j^{-\xi-1})$ and (7.1) holds with $k = 0$ for $\xi > 0$,

$$\sum_{j=0}^s \Delta_j(-\xi)\{\log(t-j)\}^r (t-j)^\gamma \sim (\log t)^r t^\gamma \sum_{j=0}^s \Delta_j(-\xi) = O(t^{\gamma+\delta_1}s^{-\xi})$$

for $s = o(t)$, $\delta_1 > 0$. On the other hand,

$$\left| \sum_{j=s+1}^{t-1} \Delta_j(-\xi)\{\log(t-j)\}^r (t-j)^\gamma \right| \le Cs^{-\xi-1}(\log t)^r \sum_{j=1}^{t} j^\gamma.$$

The sum on the right-hand side is $O(t^{1+\gamma})$ for $\gamma > -1$, $O((\log t))$ for $\gamma = -1$ and $O(1)$ for $\gamma < -1$. Thus choosing $s = t^{1-\delta_2/(\xi+1)}$, $\delta_2 > 0$, produces the result.                                                                                    □

LEMMA 2.   *For $w_t = t^\gamma$ and any integer $r > 0$, as $t \to \infty$*

$$(7.7) \qquad (-\log \Delta)^r w_t^\# \sim (\log t)^r t^\gamma \qquad\qquad\qquad \textit{for } \gamma > -1,$$

$$(7.8) \qquad\qquad = O(t^{-1}(\log t)^{r-1}\{\mathbb{1}(\gamma < -1) + (\log t)\mathbb{1}(\gamma = -1)\})$$

$$\textit{for } \gamma \le -1.$$

PROOF.   Suppose (7.7) is true for a given $r$. Then as $t \to \infty$

$$(7.9) \quad (-\log \Delta)^{r+1} w_t^\# \sim (-\log \Delta)(\log t)^r w_t^\# = \sum_{j=1}^{t-1} j^{-1}\{\log(t-j)\}^r (t-j)^\gamma.$$

The difference between this and

$$(7.10) \qquad\qquad\qquad (\log t)^r \sum_{j=1}^{t-1} j^{-1}(t-j)^\gamma$$

is bounded by $C(\log t)^{r-1}$ times

$$\sum_{j=1}^{t-1} j^{-1}\{\log t - \log(t-j)\}(t-j)^\gamma \le \sum_{j=1}^{t-1} j^{-1}|\log(1-j/t)|(t-j)^\gamma.$$

Splitting this into sums over $j \in [1, [t/2]]$ and $j \in [[t/2]+1, t-1]$, it is seen that the first of these is bounded by

$$t^{-1} \sum_{j=1}^{t-1} (t-j)^\gamma \le Ct^\gamma,$$

since $|\log(1-x)| \le x$ for $x \in (0, \frac{1}{2})$, while the second is bounded by

$$Ct^{-1} \sum_{j=1}^{t-1} |\log(j/t)| j^\gamma \le Ct^\gamma \log t.$$

The difference between (7.10) and

$$(7.11) \qquad\qquad\qquad (\log t)^r t^\gamma \sum_{j=1}^{t-1} j^{-1}$$

is bounded by

$$C(\log t)^r t^\gamma \sum_{j=1}^{t-1} j^{-1} |(1-j/t)^\gamma - 1| \le C(\log t)^r t^\gamma.$$

Then (7.11) $\sim (\log t)^{r+1} t^\gamma$ as $t \to \infty$. For $\gamma \le -1$, we can write

$$(-\log \Delta)^r w_t^\# = \sum_{j=1}^{t-1} a_j^{(r)} (t-j)^\gamma,$$

where $a_j^{(r)} = O(\{\log(j+1)\}^{r-1} j^{-1})$. Splitting the sum as before, the first one is $O((\log t)^r t^\gamma)$ and the second is $O((\log t)^{r-1} t^{-1})$ for $\gamma < -1$ and $O((\log t)^r t^{-1})$ for $\gamma = -1$.  $\square$

In the following four lemmas $b(e^{i\lambda})$ is taken to be a function with absolutely convergent Fourier series, and $b_j = (2\pi)^{-1} \int_{-\pi}^{\pi} b(e^{i\lambda}) e^{ij\lambda} d\lambda$.

LEMMA 3.  *For $w_t = t^\gamma$,*

$$b(B) w_t^\# \sim b(1) t^\gamma \qquad as \ t \to \infty.$$

PROOF.  The left-hand side equals $t^\gamma \sum_{j=0}^{t-1} b_j + \sum_{j=0}^{t-1} b_j \{(t-j)^\gamma - t^\gamma\}$. The first term differs by $o(t^\gamma)$ from $b(1) t^\gamma$, and the second is bounded by

$$Ct^\gamma \sum_{j=0}^{t-1} |b_j| \left| 1 - \left(1 - \frac{j}{t}\right)^\gamma \right| \le Ct^{\gamma-1} \sum_{j=0}^{t} j |b_j| = o(t^\gamma)$$

from the Toeplitz lemma.  $\square$

LEMMA 4.  *For a sequence $w_t$ such that $w_t = 0$, $t \le 0$, and any integer $r$, as $\xi \to \xi_0$*

$$(7.12) \quad (\log \Delta)^r (\Delta^\xi - \Delta^{\xi_0}) b(B) w_t = (\log \Delta)^{r+1} \Delta^{\xi_0} b(B) w_t (\xi - \xi_0)$$

$$+ O\left( \left\{ \sum_{j=1}^{t} (\Delta^m w_{t-j})^2 \right\}^{1/2} (\xi - \xi_0)^2 \right)$$

*for $m \in (\xi_0 - \frac{1}{2}, \xi_0 + \frac{1}{2})$.*

PROOF.  By the mean value theorem the left-hand side of (7.12) is

$$(\log \Delta)^{r+1} \Delta^{\xi_0} b(B) w_t (\xi - \xi_0) + \tfrac{1}{2} (\log \Delta)^{r+2} b(B) \Delta^{\bar{\xi}} w_t (\xi - \xi_0)^2,$$

for $|\bar{\xi} - \xi_0| \le |\xi - \xi_0|$. The last term can be written $\tfrac{1}{2} \sum_{j=1}^{t-1} c_j \Delta^m w_{t-j} (\xi - \xi_0)^2$, where $c_j$ is the coefficient of $s^j$ in the Taylor expansion of $\{\log(1-s)\}^{r+2} \times$ $(1-s)^{\bar\xi - m}$. From Stirling's approximation, $c_j \sim (\log j)^{r+2} j^{m-\bar\xi-1}$ as $j \to \infty$. Now $m - \bar\xi \le m - \xi_0 + |\xi - \xi_0|$. The right-hand side of this is less than $\frac{1}{2}$ if $|\xi - \xi_0| < \frac{1}{2} - m + \xi_0$, where the right-hand side of the latter inequality is positive. Thus for $|\xi - \xi_0|$ small enough, $m - \bar\xi - 1 < -\frac{1}{2}$. Then $\sum_{j=1}^\infty c_j^2 < \infty$ for all $r$, so the proof is completed by the Cauchy inequality. $\square$

LEMMA 5.   *For real $\xi$ and $m_0$ defined by (1.2),*

$$(7.13) \qquad \Delta^\xi b(B) x_t = \Delta^{\xi - m_0} b(B) v_t^\#, \qquad t \in \mathbb{Z}.$$

PROOF.   The left-hand side of (7.13) is

$$\Delta^\xi b(B) \Delta^{-m_0} v_t^\# = \Delta^{\xi - m_0} b(B) v_t^\#, \qquad t \in \mathbb{Z}. \qquad \square$$

The next lemma gives a uniform bound for the variance of a process that is only "asymptotically stationary."

LEMMA 6.   *For all $r \ge 0$, and $\zeta_0$ defined by (1.4),*

$$(7.14) \qquad E\{(-\log \Delta)^r \Delta^{\zeta_0} b(L) v_t^\#\}^2 \le C < \infty.$$

PROOF.   The left-hand side of (7.14) is

$$(7.15) \qquad \int_{-\pi}^\pi \left| \sum_{j=0}^{t-1} c_j e^{ij\lambda} \right|^2 |1 - e^{i\lambda}|^{-2\zeta_0} f(\lambda)\, d\lambda \le C \left( \sum_{j=0}^\infty |c_j| \right)^2$$

for $\zeta_0 > 0$ since $|1 - e^{i\lambda}|^{-2\zeta_0} f(\lambda)$ is integrable, $c_j$ being the $j$th Fourier coefficient of $[\{-\log(1 - e^{i\lambda})\}^r (1 - e^{i\lambda})^{\zeta_0}] b(e^{i\lambda})$. The $j$th Fourier coefficient of the factor in braces is $O((\log j)^r j^{-\zeta_0 - 1})$, so since the $b_j$ are summable so are the $c_j$. For $\zeta_0 \le 0$ $|1 - e^{i\lambda}|^{-2\zeta_0} f(\lambda)$ is bounded, so the left-hand side of (7.15) is bounded by $\sum_0^\infty c_j^2 < \infty$. $\square$

LEMMA 7.   *Let $S_m$ be the $m \times m$ matrix with $(j,k)$th element $(j, k \ge 1)$,*

$$\int_{-1}^1 u^{j+k-2}\, du = 2(j+k-1)^{-1} \mathbb{1}(j+k \text{ even}).$$

*Then for $m$ sufficiently large,*

$$\mathrm{tr}(S_m^{-1}) < (2\pi)^{-2} \left[ \frac{8}{3} + \frac{1}{2} \log\left\{ (2m-3)\left(\frac{2m}{3} - 1\right)\right\}\right] \eta^{2m}.$$

PROOF.   It is clear that, like $S_m$, $S_m^{-1}$ must have $(j, k)$th element that is zero for all odd $j + k$. This immediately ensures the necessary property that even rows (columns) of $S_m$ are orthogonal to odd rows (columns) of $S_m^{-1}$. It then suffices to study the two square matrices $S_{1,m}$ and $S_{2,m}$ formed from, respectively, the odd and even rows and columns of $S_m$. These exclude only and all zero elements of $S_m$, and $S_m^{-1}$ is the $m \times m$ matrix whose $(2j-1, 2k-1)$th element is the $(j+k)$th element of $S_{1,m}^{-1}$, whose $(2j, 2k)$th element is the $(j,k)$th element of $S_{2,m}^{-1}$, and whose other elements are all zero. Thus it suffices to consider $S_{1,m}^{-1}$ and $S_{2,m}^{-1}$, and indeed $\text{tr}(S_m^{-1}) = \text{tr}(S_{1,m}^{-1}) + \text{tr}(S_{2,m}^{-1})$. We take $m$ to be even; details for $m$ odd are only slightly different and since we want a result only for large $m$ this outcome will clearly be unaffected.

$S_{1,m}$ and $S_{2,m}$ are both Cauchy matrices (see, e.g., [17], page 36), having $(j,k)$th element of the form $(a_j + a_k)^{-1}$, in particular, $(j + k - \frac{3}{2})^{-1}$, $(j + k - \frac{1}{2})^{-1}$, respectively. From Knuth [17], page 36, the $j$th diagonal elements of $S_{1,m}^{-1}$, $S_{2,m}^{-1}$ are, respectively, $2U_1^2(j)/(4j-3)$, $2U_2^2(j)/(4j-1)$, where we define, for real $s$,

$$U_1(s) = \frac{\prod_{1 \leq i \leq m/2}(i + s - 3/2)^2}{\prod_{1 \leq i \leq m/2, i \neq s}(i - s)},$$

$$U_2(s) = \frac{\prod_{1 \leq i \leq m/2}(i + s - 1/2)^2}{\prod_{1 \leq i \leq m/2, i \neq j}(i - s)}.$$

Thus

$$\begin{aligned}
\text{tr}(S_m^{-1}) &= 2 \sum_{j=1}^{m/2} \{(4j-3)^{-1}U_1^2(j) + (4j-1)^{-1}U_2^2(j)\} \\
&\leq \left\{2 + \frac{1}{2}\log(2m-3)\right\} \max_{1 \leq j \leq m/2} U_1^2(j) \\
&\quad + \left\{\frac{2}{3} + \frac{1}{2}\log\left(\frac{2m}{3} - \frac{1}{3}\right)\right\} \max_{1 \leq j \leq m/2} U_2^2(j).
\end{aligned}$$

For $s \in (0, m/2 - 1)$

$$U_1(s) - U_1(s+1) = U_1(s)\left\{1 - \frac{(s + m/2 - 1/2)(m/2 - s)}{(s - 1/2)s}\right\}.$$

The factor in braces is $2 - m(m-1)/\{2s(2s-1)\}$, which is negative for $s < s(m)$ and positive for $s > s(m)$, where $s(m) = \frac{1}{4} + \{2m(m-1) + 1\}^{1/2}/4 \sim m/\sqrt{8}$ as $m \to \infty$. Thus, as $m \to \infty$

$$(7.16) \quad \max_{1 \leq j \leq m/2} U_1(j) \sim \frac{\Gamma((1/2 + 1/\sqrt{8})m - 1/2)}{\Gamma(m/\sqrt{8} - 1/2)\Gamma(m/\sqrt{8})\Gamma((1/2 - 1/\sqrt{8})m + 1)}.$$

Applying Stirling's approximation, that is,

$$\Gamma(am + b) \sim (2\pi)^{1/2} e^{-am}(am)^{am+b-1/2}$$

as $m \to \infty$, and noting that

$$\left\{ \frac{(1+2^{-1/2})^{1+2^{-1/2}} 2^{2^{-1/2}}}{(1-2^{-1/2})^{1-2^{-1/2}}} \right\}^{1/2} = 1 + 2^{1/2},$$

(7.16) is $(2\pi)^{-1}\eta^m(1+o(1))$. In the same way it can be seen that $U_2(s)$ is maximized at $\{2m(m+1)+1\}^{1/2}/4 - \frac{1}{4} \sim m/\sqrt{8}$, whence $\max_{1 \le j \le m/2} U_2(j) \sim (2\pi)^{-1}\eta^m(1+o(1))$ also. The proof is then routinely completed.   □

Denote by $\underline{\lambda}(A)$ the smallest eigenvalue of the matrix $A$.

LEMMA 8.   *As $L \to \infty$,*

$$\underline{\lambda}(W^{(L)})^{-1} = O(\pi_L).$$

PROOF.   The method of proof, given Lemma 7, is similar to one in [25], but we obtain a refinement. Define $\phi_+^{(L)}(s) = (1, \phi^{(L)}(s)^T)^T$, $W_+^{(L)} = E\{\phi_+^{(L)}(\varepsilon_t) \times \phi_+^{(L)}(\varepsilon_t)^T\}$, so $W^{(L)} = PW_+^{(L)}P^T$, where the $L \times (L+1)$ matrix $P$ consists of the last $L$ rows of the $(L+1)$-rowed identity matrix. Then $\underline{\lambda}(W^{(L)}) \ge \underline{\lambda}(W_+^{(L)})\underline{\lambda}(PP^T) = \underline{\lambda}(W_+^{(L)})$. If $(-1,1) \subset (\phi(s_1), \phi(s_2))$ (which implies $\varphi \le 1$), then [since $\phi'(s)$ is bounded on $(s_1, s_2)$] $\underline{\lambda}(W_+^{(L)}) \ge \underline{\lambda}(S_{L+1})/C \ge \mathrm{tr}(S_{L+1}^{-1})^{-1}/C$, where we use $S_m$ defined as in Lemma 7, which can then be applied. Otherwise, $W_+^{(L)}$ exceeds, by a nonnegative definite matrix,

$$(7.17) \quad C^{-1} \int_{\phi(s_1)}^{\phi(s_2)} u^{(L)} u^{(L)T} \, du = \left\{ \frac{\phi(s_2) - \phi(s_1)}{C} \right\} A \int_{-1}^{1} u^{(L)} u^{(L)T} \, du \, A^T,$$

where $u^{(L)} = (1, u, \dots, u^L)^T$ and $A$ is the lower-triangular matrix with $(i,j)$th element $\binom{i-1}{j-i} \phi(s_1)^{i-j} \{\phi(s_2) - \phi(s_1)\}^{j-1}$, $j \le i$. The smallest eigenvalue of (7.17) is no less than $C^{-1}\{\phi(s_2) - \phi(s_1)\}\underline{\lambda}(AA^T)\underline{\lambda}(S_{L+1})$. Now $\underline{\lambda}(AA^T) \ge \|A^{-1}\|^{-2}$, where by recursive calculation $A^{-1}$ is seen to be lower-triangular with $(i,j)$th element $a^{ij} = \binom{i-1}{j-i}\{-\phi(s_1)\}^{i-j}\{\phi(s_2) - \phi(s_1)\}^{1-i}$, $j \le i$. Thus

$$\|A^{-1}\|^2 = \sum_{i=1}^{L+1} \left( \sum_{j=1}^{i} a^{ij2} \right) \le \sum_{i=1}^{L+1} \left( \sum_{j=1}^{i} |a^{ij}| \right)^2 \le \sum_{i=1}^{L+1} \varphi^{2(i-1)}.$$

This is bounded by $(1-\varphi^2)^{-1}$ for $\varphi < 1$, by $L+1$ for $\varphi = 1$ and by $(\varphi^2 - 1)^{-1} \times \varphi^{2(L+1)}$ for $\varphi > 1$.   □

LEMMA 9.   *For $a \ge 0, b \ge 0$,*

$$(7.18) \quad \sum_{\ell=1}^{L} \mu_{a\ell+b} \le \rho_{aL}.$$

PROOF. In case $a = 0$, or $a > 0$ but Assumption A2(c) holds, this is trivial. For $a > 0$ under Assumption A2(b), monotonic nondecrease of $\mu_a$ in real $a$ implies that the left-hand side of (7.18) is bounded by

$$C \sum_{\ell=1}^{[aL+b]} \mu_{\kappa\ell} \leq \left( \frac{CL}{t} \right)^{(a/\kappa)L} E(e^{t|\varepsilon_0|^\kappa})$$

for any $t \in (0, 1)$, and by Assumption A2(b) there exists such $t$ that this is bounded by $\rho_{aL}$. $\square$

LEMMA 10. As $n \to \infty$,

$$\begin{aligned}
(7.19) \qquad & \|a^{(L)}\| = O(L\rho_{2\kappa L}^{1/2} \pi_L), \\
& \|\hat{a}^{(L)}(\varepsilon) - a^{(L)}\| = O\left( \frac{L}{n^{1/2}} \rho_{2\kappa L}^{1/2} \pi_L (1 + L^{1/2} \rho_{4\kappa L}^{1/2} \pi_L) \right).
\end{aligned}$$

PROOF. Write

$$\hat{a}^{(L)}(\varepsilon) - a^{(L)} = \{ W^{(L)}(\varepsilon)^{-1} - W^{(L)-1} \} w^{(L)}(\varepsilon) + W^{(L)-1} \{ w^{(L)}(\varepsilon) - w^{(L)} \}.$$

From (6.23), the Schwarz inequality and Lemma 9

$$\|w^{(L)}\|^2 = \sum_{\ell=1}^{L} \ell^2 \{ E\{ \phi'(\varepsilon_0) \phi^{\ell-1}(\varepsilon_0) \} \}^2 \leq CL^2 \sum_{\ell=1}^{L} \mu_{2\kappa(\ell+K)} \leq L^2 \rho_{2\kappa L}.$$

Similarly, and from independence of the $\varepsilon_t$,

$$E\|w^{(L)}(\varepsilon) - w^{(L)}\|^2 \leq n^{-1} \sum_{\ell=1}^{L} \ell^2 E\{ \phi'(\varepsilon_0) \phi^{\ell-1}(\varepsilon_0) \}^2 \leq (L^2/n) \rho_{2\kappa L},$$

$$E\|W^{(L)}(\varepsilon) - W^{(L)}\|^2 \leq n^{-1} \sum_{k,\ell=1}^{L} E\{ \phi(\varepsilon_0)^{2(k+\ell)} \} \leq (L/n) \rho_{4\kappa L}.$$

Now apply Lemma 8. $\square$

LEMMA 11. For $j \geq 0$ let $\alpha_j = \Delta_j(d)$ for $d \leq 1$ and $|\beta_j| \leq C(j+1)^{-3}$. Then the sequence $\sum_{k=0}^{j} \alpha_{j-k} \beta_k$, $j \geq 0$, has property $P_0(d)$.

PROOF. By Stirling's approximation $\alpha_j$ has property $P_0(d)$, whence the proof is completed by splitting sums around $j/2$ and elementary bounding of each. $\square$

LEMMA 12.   *For $j \geq 0$ let the sequence $\alpha_j$, $j \geq 0$, have property $P_0(-d)$ and for $d > 0$ let $\sum_{j=0}^{\infty} \alpha_j = 0$. Then for $|d| < 1$ the sequence*

$$\gamma_j = \sum_{k=0}^{j} (j+1-k)^{-1} \alpha_k, \qquad j \geq 0,$$

*has property $P_1(-d)$.*

PROOF.   We give the proof only of $|\gamma_j - \gamma_{j+1}| \leq C\{\log(j+1)\}j^{-d-2}$, the proof of $|\gamma_j| \leq C\{\log(j+1)\}j^{-d-1}$ being similar and simpler. We have

$$\gamma_j - \gamma_{j+1} = \sum_{k=0}^{\tilde{j}} \{(j+1-k)^{-1} - (j+2-k)^{-1}\} \alpha_k - (j+1-\tilde{j})^{-1} \alpha_{\tilde{j}+1}$$

$$+ \sum_{k=\tilde{j}+1}^{j} (j+1-k)^{-1} (\alpha_k - \alpha_{k+1}),$$

where $\tilde{j} = [j/2]$. The second term is bounded by $Cj^{-d-2}$ and the third by $C(\log j)j^{-d-2}$. For $d < 0$ the first term is bounded by $Cj^{-d-2}$ and for $d = 0$ by $C(\log j)j^{-d-2}$. For $d > 0$ we apply summation by parts to this first term and $\sum_{j=0}^{\infty} \alpha_j = 0$ to obtain the bound $Cj^{-d-2}$ again.   □

LEMMA 13.   *Let the sequence $\alpha_j$, $j \geq 0$, have property $P_0(-d)$ and the sequence $\beta_j$, $j \geq 0$, have property $P_0(e)$, and let*

$$(7.20) \qquad \begin{aligned} \sum_{j=0}^{\infty} |\alpha_j| &< \infty \qquad \text{if } d = 0, \\ \sum_{j=0}^{\infty} |\beta_j| &< \infty \qquad \text{if } e = 0, \\ \sum_{j=0}^{\infty} \beta_j &= 0 \qquad \text{if } e < 0. \end{aligned}$$

*Then for $|d| < 1$, $|e| < 1$ it follows that for all $j > 0$, $t > 0$,*

$$(7.21) \qquad \left| \sum_{k=0}^{j} \alpha_{k+t} \beta_{j-k} \right| \leq Cj^e t^{-d-1}, \qquad\qquad j \leq t,$$

$$(7.22) \qquad\qquad\qquad\qquad \leq Cj^{e-1} \max(j^{-d}, t^{-d}), \qquad j > t.$$

*If instead $\alpha_j$ has property $P_1(-d)$ and (7.20) is not imposed,*

$$(7.23) \qquad \left| \sum_{k=0}^{j} \alpha_{k+t} \beta_{j-k} \right| \leq C(\log^{r+1} t) j^e t^{-d-1}, \qquad\qquad j \leq t,$$

$$(7.24) \qquad \leq C(\log^{r+1} j) j^{e-1} \max(j^{-d}, t^{-d}), \qquad j > t.$$

PROOF. We prove only (7.21) and (7.22), the proof of (7.23) and (7.24) being very similar but notationally slightly more complex and less elegant. Write $S_{ab} = \sum_{k=a}^{b} \alpha_{t+k} \beta_{j-k}$. We have

$$|S_{0j}| \leq t^{-d-1} \sum_{k=0}^{j} |\beta_k| \leq C j^e t^{-d-1}, \qquad e \geq 0.$$

This proves (7.21) for $e \geq 0$ and all $d$. On the other hand, with $\tilde{j} = [j/2]$, summation by parts gives

$$|S_{0\tilde{j}}| \leq \sum_{k=0}^{\tilde{j}-1} |\beta_{j-k} - \beta_{j-k-1}| \sum_{i=0}^{k} |\alpha_{t+i}| + |\beta_{j-\tilde{j}}| \sum_{k=0}^{\tilde{j}} |\alpha_{t+k}|$$

$$(7.25) \qquad \leq C t^{-d} \left\{ \sum_{k=0}^{\tilde{j}} (j-k)^{e-2} + j^{e-1} \right\}$$

$$\leq C j^{e-1} t^{-d}, \qquad d \geq 0, \text{all } e,$$

while

$$(7.26) \qquad |S_{\tilde{j}+1,j}| \leq C(t+\tilde{j})^{-d-1} j^e \leq C j^{e-d-1} \qquad \text{all } d; e \geq 0.$$

This proves (7.22) for $d \geq 0$, $e \geq 0$ since $j^{e-d-1} \leq j^{e-1} t^{-d}$, $j > t$. For $e < 0$

$$S_{0j} = -\sum_{k=0}^{j-1} \{\alpha_{j-k+t} - \alpha_{j-k-1+t}\} \sum_{i=k+1}^{\infty} \beta_i - \alpha_t \sum_{k=j+1}^{\infty} \beta_i,$$

since $\sum_{j=0}^{\infty} \beta_j = 0$. This is bounded by $C\{t^{-d-2} j^{e+1} + t^{-d-1} j^e\} \leq C j^e t^{-d-1}$ for $j \leq t$, to prove (7.21) for $e < 0$ and all $d$. For $e < 0$ and all $d$

$$S_{\tilde{j}+1,j} = \sum_{k=0}^{j-\tilde{j}-1} \alpha_{j+t-k} \beta_k$$

$$= -\sum_{k=0}^{j-\tilde{j}-2} (\alpha_{j+t-k} - \alpha_{j+t-k-1}) \sum_{i=k+1}^{\infty} \beta_i - \alpha_{t+\tilde{j}-1} \sum_{k=j-\tilde{j}}^{\infty} \beta_k,$$

and this is bounded by $C\{(t+\tilde{j})^{-d-2} j^{e+1} + (t+\tilde{j})^{-d-1} j^e\} \leq C j^{e-1} t^{-d}$, which with (7.25) proves (7.22) for $d \geq 0$, $e < 0$. Finally, for $d < 0$ and all $e$

$$|S_{0\tilde{j}}| = \left| \sum_{k=j-\tilde{j}}^{j} \alpha_{j+t-k} \beta_k \right| \leq C j^{e-d-1},$$

which with (7.26) completes the proof of (7.22). $\square$

LEMMA 14.    *For $|\zeta_0| < \frac{1}{2}$,*

$$\sum_{j=0}^{\infty} \lambda_{jt}^2 \le Ct^{-1}, \tag{7.27}$$

$$\sum_{j=0}^{\infty} \left( \sum_{t=1}^{n} \lambda_{jt} \right)^2 \le Cn. \tag{7.28}$$

PROOF.    In this and subsequent proofs we drop the zero subscript from $\zeta_0$. We omit the proof for $\zeta = 0$ as it is simple. From Lemma 13

$$\sum_{j=1}^{\infty} \lambda_{jt}^2 \le Ct^{-2\zeta-2} \sum_{j=1}^{t} j^{2\zeta} + C \sum_{j=t}^{\infty} j^{2\zeta-2} \max(j^{-2\zeta}, t^{-2\zeta}).$$

The first sum is bounded by $Ct^{2\zeta+1}$ and the second by $Ct^{-2\zeta} \sum_{j=t}^{\infty} j^{2\zeta-2} \le Ct^{-1}$ when $\zeta > 0$ and by $C \sum_{j=t}^{\infty} j^{-2} \le Ct^{-1}$ when $\zeta < 0$, to prove (7.27). For $j < n$ and $\zeta \ne 0$

$$\left| \sum_t \lambda_{jt} \right| \le Cj^{\zeta-1} \sum_{t=1}^{j} \max(j^{-\zeta}, t^{-\zeta}) + Cj^{\zeta} \sum_{t=j+1}^{n} t^{-\zeta-1}$$

$$\le C \max(1, (j/n)^{\zeta}).$$

For $j \ge n$

$$\left| \sum_t \lambda_{jt} \right| \le Cj^{\zeta-1} \sum_{t=1}^{n} \max(j^{-\zeta}, t^{-\zeta}) \le C \max(n/j, (n/j)^{1-\zeta}).$$

Thus

$$\sum_{j=0}^{\infty} \left( \sum_t \lambda_{jt} \right)^2 = \sum_{j=0}^{n} \left( \sum_t \lambda_{jt} \right)^2 + \sum_{j=n+1}^{\infty} \left( \sum_t \lambda_{jt} \right)^2$$

$$\le Cn + Cn^{2-2\zeta} \sum_{j=n}^{\infty} j^{2\zeta-2} \le Cn, \qquad \zeta > 0,$$

$$\le Cn^{-2\zeta} \sum_{j=1}^{n} j^{2\zeta} + n^2 \sum_{j=n}^{\infty} j^{-2} \le Cn, \qquad \zeta < 0,$$

to prove (7.28).    □

Define

$$h_{jk} = \sum_t (t+j)^{-1} |\lambda_{kt}|, \qquad j, k \ge 1.$$

LEMMA 15. *For* $0 < \zeta_0 < \frac{1}{2}$ *and* $j \geq 1$,

$$(7.29) \quad h_{jk} \leq Cj^{-1/2} \min(j^{-1/2}, k^{-1/2}), \qquad 1 \leq k \leq n,$$

$$(7.30) \qquad \leq Cj^{-1} k^{\zeta_0 - 1} n^{1/2 - \zeta_0} \min(j^{1/2}, n^{1/2}), \qquad k \geq n.$$

*For* $-\frac{1}{2} < \zeta_0 \leq 0$ *and* $j \geq 1$,

$$h_{jk} \leq C \min(j^{1/2 - \varepsilon} k^{-1/2 + \varepsilon}, k^{-1} \log k),$$

$$(7.31)$$

$$0 < \varepsilon < \frac{1}{2} + \zeta_0, 1 \leq k < n,$$

$$(7.32) \qquad \leq Ck^{-1} \min(n/j, \log n), \qquad k \geq n.$$

PROOF. It follows from Lemma 13 that for $1 \leq k \leq n$,

$$(7.33) \quad h_{jk} \leq Ck^{\zeta - 1} \sum_{t=1}^{k} (t+j)^{-1} \max(k^{-\zeta}, t^{-\zeta}) + Ck^{\zeta} \sum_{t=k}^{n} (t+j)^{-1} t^{-\zeta - 1}.$$

Suppose $\zeta > 0$. The first term on the right-hand side is bounded by

$$Cj^{-1} k^{\zeta - 1} \sum_{t=1}^{k} t^{-\zeta} \leq Cj^{-1}, \qquad j \geq k,$$

$$Cj^{-1/2} k^{\zeta - 1} \sum_{t=1}^{k} t^{-\zeta - 1/2} \leq C(jk)^{-1/2}, \qquad j \leq k.$$

The second term on the right-hand side of (7.33) is bounded by

$$Cj^{-1} k^{\zeta} \sum_{t=k}^{n} t^{-\zeta - 1} \leq Cj^{-1}, \qquad j \geq k,$$

$$Cj^{-1/2} k^{\zeta} \sum_{t=k}^{n} t^{-\zeta - 3/2} \leq C(jk)^{-1/2}, \qquad j \leq k.$$

This proves (7.29). Let $\zeta \leq 0$. The first term on the right-hand side of (7.33) is bounded by

$$Ck^{-1} \sum_{t=1}^{k} (t+j)^{-1} \leq C \min(j^{-1}, k^{-1} \log k)$$

and the second by

$$Ck^{\zeta} j^{-1/2 - \varepsilon} \sum_{t=k}^{\infty} t^{-\zeta - 3/2 + \varepsilon} \leq Cj^{-1/2 - \varepsilon} k^{-1/2 + \varepsilon}, \qquad j \geq k,$$

$$Ck^{\zeta} \sum_{t=k}^{n} t^{-\zeta - 2} \leq Ck^{-1}, \qquad j \leq k.$$

This proves (7.31). For $k \geq n$ (7.30) and (7.32) are readily deduced from

$$h_{jk} \leq Ck^{\zeta-1} \sum_t (t+j)^{-1} t^{-\zeta} \mathbb{1}(\zeta > 0) + Ck^{-1} \sum_t (t+j)^{-1} \mathbb{1}(\zeta \leq 0). \qquad \square$$

LEMMA 16. *For* $|\zeta_0| < \frac{1}{2}$,

$$E \left\| \sum_t \varepsilon'_{t1} \sum_{j=0}^{\infty} \lambda_{jt} \varepsilon_{-j} \right\|^2 \leq C(\log n)^3.$$

PROOF. Writing $\gamma(s; \nu_0) = \sum_{j=0}^{\infty} \gamma_j s^j$, the expression within the norm is

$$(7.34) \qquad \sum_t \sum_{j=1-t}^{-1} \gamma_{t+j} \varepsilon_{-j} \sum_{k=0}^{\infty} \lambda_{kt} \varepsilon_{-k} + \sum_{j,k=0}^{\infty} H_{jk} \varepsilon_{-j} \varepsilon_{-k},$$

where $H_{jk} = \sum_t \gamma_{j+t} \lambda_{kt}$. The squared norm of the first term has expectation bounded by

$$\sum_s \sum_t \left( \sum_{j=\max(1-s,1-t)}^{-1} \|\gamma_{s+j}\| \|\gamma_{t+j}\| \right) \left( \sum_{k=0}^{\infty} \lambda_{sk} \lambda_{tk} \right).$$

For $s \leq t$ the first bracketed factor is $O((t-s+1)^{-1} \log n)$ because $\|\gamma_j\| \leq C(j+1)^{-1}$, while the second one is bounded by

$$Ct^{-\zeta-1} s^{-\zeta-1} \sum_{j=1}^{s} j^{2\zeta} + Ct^{-\zeta-1} \sum_{j=s+1}^{t} j^{2\zeta-1} \max(j^{-\zeta}, s^{-\zeta})$$

$$+ C \sum_{j=t+1}^{\infty} j^{2\zeta-2} \max(j^{-\zeta}, s^{-\zeta}) \max(j^{-\zeta}, t^{-\zeta})$$

$$\leq C\{s^{-\zeta} t^{\zeta-1} \mathbb{1}(\zeta > 0) + s^{\zeta} t^{-\zeta-1} \mathbb{1}(\zeta < 0) + (st)^{-1/2} \mathbb{1}(\zeta = 0)\}$$

$$\leq C(st)^{-1/2}.$$

We have

$$\sum_{s=1}^{t} (t-s+1)^{-1} s^{-1/2} \leq \sum_{s=1}^{[t/2]} (t-s+1)^{-1} s^{-1/2} + \sum_{s=[t/2]}^{t} (t-s+1)^{-1} s^{-1/2}$$

$$\leq C(\log t) t^{-1/2},$$

$$C(\log n) \sum_t (\log t) t^{-1} \leq C(\log n)^3.$$

Next, since $|H_{jk}| \leq Ch_{jk}$, the squared norm of the second term on the right-hand side of (7.34) has expectation bounded by

$$C \sum_{j,k=0}^{\infty} (h_{jk}^2 + h_{jj} h_{kk} + h_{jk} h_{kj}).$$

We apply Lemma 15 to complete the proof. For $\zeta > 0$

$$\sum_{j,k=0}^{\infty} h_{jk}^2 \leq C \sum_{k=1}^{n} \sum_{j=1}^{k} (jk)^{-1} + C \sum_{k=1}^{n} \sum_{j=k}^{\infty} j^{-2}$$

$$+ Cn^{1-2\zeta} \sum_{k=n}^{\infty} \sum_{j=1}^{n} j^{-1} k^{2\zeta-2} + Cn^{2-2\zeta} \sum_{k=n}^{\infty} \sum_{j=n}^{\infty} j^{-2} k^{2\zeta-2}$$

$$\leq C(\log n)^2,$$

$$\sum_{j=0}^{\infty} h_{jj} \leq C \sum_{j=1}^{n} j^{-1} + n^{1-\zeta} \sum_{j=n}^{\infty} j^{-2} \leq C \log n$$

and

$$\sum_{j,k=0}^{\infty} h_{jk} h_{kj} \leq C \sum_{k=1}^{n} \sum_{j=1}^{k} (jk)^{-1} + Cn^{1/2-\zeta} \sum_{k=1}^{n} \sum_{j=k}^{\infty} j^{\zeta-2} k^{-1/2}$$

$$+ Cn^{2-2\zeta} \sum_{j,k=n}^{\infty} (jk)^{\zeta-2}$$

$$\leq C(\log n)^2.$$

For $\zeta \leq 0$

$$\sum_{j,k=0}^{\infty} h_{jk}^2 \leq \sum_{k=1}^{n} \sum_{j=1}^{k} (k^{-1} \log k)^2 + C \sum_{k=1}^{n} \sum_{j=k}^{\infty} j^{-1-2\varepsilon} k^{-1+2\varepsilon}$$

$$+ C(\log n)^2 \sum_{k=n}^{\infty} \sum_{j=1}^{n} k^{-2} + Cn^2 \sum_{j,k=n}^{\infty} (jk)^{-2}$$

$$\leq C(\log n)^3,$$

$$\sum_{j=0}^{\infty} h_{jj} \leq C \sum_{j=1}^{n} j^{-1} + Cn \sum_{j=n}^{\infty} j^{-2} \leq C \log n,$$

$$\sum_{j,k=0}^{\infty} h_{jk} h_{kj} \leq C \sum_{k=1}^{n} \sum_{j=1}^{k} j^{-1/2+\varepsilon} k^{-3/2-\varepsilon} \log k$$

$$+ C \log n \sum_{k=1}^{n} \sum_{j=n}^{\infty} j^{-1/2-\varepsilon} k^{-1/2+\varepsilon} j^{-1} + Cn^2 \sum_{j,k=n}^{\infty} (jk)^{-2}$$

$$\leq C(\log n)^2. \qquad \square$$

LEMMA 17.

$$E\left\|\sum_t \varepsilon'_{t1}\right\|^4 \leq C(\log n)^4 n^2.$$

PROOF. We have

$$\sum_t \varepsilon'_{t1} = \sum_{j=1}^{n-1}\left(\sum_{i=1}^{n-j}\gamma_i\right)\varepsilon_j + \sum_{j=0}^{\infty}\left(\sum_{i=j+1}^{j+n}\gamma_i\right)\varepsilon_{-j}.$$

Thus

$$E\left\|\sum_t \varepsilon'_{t1}\right\|^4 \leq C\left(\sum_{j=1}^{n-1}\left\|\sum_{i=1}^{n-j}\gamma_i\right\|^2\right)^2 + C\left(\sum_{j=0}^{\infty}\left\|\sum_{i=j+1}^{j+n}\gamma_i\right\|^2\right)^2.$$

Since

$$\left\|\sum_{i=1}^{n-j}\gamma_i\right\| \leq \sum_{i=1}^{n-j}\|\gamma_i\| \leq C\sum_{i=1}^{n}i^{-1} \leq C\log n, \qquad 1 \leq j < n,$$

$$\left\|\sum_{i=j+1}^{j+n}\gamma_i\right\| \leq C\sum_{i=j+1}^{j+n}i^{-1} \leq C\log n, \qquad 1 \leq j \leq n,$$

$$\leq Cn/j, \qquad\qquad\qquad j \geq n,$$

the proof is readily completed.   □

LEMMA 18.  *For any sequence* $c_j$, $j \geq 0$, *and any* $r \geq 1$, *if* $\mu_{r+} < \infty$,

$$E\left|\sum_{j=0}^{\infty}c_j\varepsilon_{-j}\right|^r \leq (Cr)^{2r}\left(\sum_{j=0}^{\infty}c_j^2\right)^{r/2}\mu_{r+}^{r/r_+},$$

*where* $r_+$ *is the smallest even integer such that* $r_+ \geq r$.

PROOF.  For $r \leq 2$ the proof follows by Jensen's inequality and direct calculation. For $r > 2$ the Marcinkiewicz–Zygmund inequality indicates that

$$(7.35) \qquad E\left|\sum_{j=0}^{\infty}c_j\varepsilon_{-j}\right|^r \leq C_r E\left(\sum_{j=0}^{\infty}c_j^2\varepsilon_{-j}^2\right)^{r/2},$$

where $C_r = \{18r^{3/2}(r-1)^{-1/2}\}^r$ (see [13], page 23). By the $c_r$-inequality (7.35) is bounded by

$$C_r 2^{r/2-1}\left\{E\left|\sum_{j=0}^{\infty}c_j^2(\varepsilon_{-j}^2-1)\right|^{r/2} + \left(\sum_{j=0}^{\infty}c_j^2\right)^{r/2}\right\}$$

$$\leq C_r 2^{r/2-1}\left\{C_{r/2}E\left|\sum_{j=0}^{\infty}c_j^4(\varepsilon_{-j}^2-1)^2\right|^{r/4} + \left(\sum_{j=0}^{\infty}c_j^2\right)^{r/2}\right\}.$$

For $2 < r \leq 4$ the first expectation in the last line is bounded by

$$\left\{ E \sum_{j=0}^{\infty} c_j^4 (\varepsilon_{-j}^2 - 1)^2 \right\}^{r/4} \leq \left( \sum_{j=0}^{\infty} c_j^4 E \varepsilon_0^4 \right)^{r/4} \leq \left( \sum_{j=0}^{\infty} c_j^2 \right)^{r/2} \mu_4^{r/4}.$$

For $r > 4$ we instead apply the $c_r$-inequality to that expectation, and then the Marcinkiewicz–Zygmund inequality again, and so on, eventually bounding (7.35) by

$$C_r C_{r/2} C_{r/4} \cdots C_2 \cdot 2^{r/2} \cdot 2^{r/4} \cdot 2^{r/8} \cdots 1 \left( \sum_{j=0}^{\infty} c_j^2 \right)^{r/2} \mu_{r_+}^{r/r_+}.$$

The result follows on noting that $r \cdot r^{1/2} \cdot r^{1/4} \cdots r^{1/r} < r^2$, $2^{1/2} \cdot 2^{1/4} \cdots 1 < 2$, $2^{1/2} \cdot 4^{1/4} \cdots r^{1/r} > 1$ and $j/(j-1) \leq 2$ for all $j \geq 2$. $\quad\square$

LEMMA 19. *As $n \to \infty$*

$$\|\hat{a}^{(L)}(E/\sigma_0) - \hat{a}^{(L)}(\varepsilon)\| = O_p(\rho_{2\kappa L}^{3/2} \pi_L^2 (L^2 n^{-1/2} + (CL)^{4\kappa L+3} n^{-1} \log n)).$$

PROOF. Because the proof is similar to details in Section 3 we sketch it. It turns out that $\{W^{(L)}(E/\sigma_0)^{-1} - W^{(L)}(\varepsilon)^{-1}\} w^{(L)}(E/\sigma_0)$ dominates $W^{(L)}(\varepsilon)^{-1} \times \{w^{(L)}(E/\sigma_0) - w^{(L)}(\varepsilon)\}$, so we look only at the former. $\|W^{(L)}(E/\sigma_0) - W^{(L)}(\varepsilon)\|$ is bounded by

$$(7.36) \qquad Cn^{-1} \left[ \sum_{k,\ell=1}^{L} \left\{ \left( \sum_t \delta_{kt} \delta_{\ell t} \right)^2 + \left( \sum_t \phi_k(\varepsilon_t) \delta_{\ell t} \right)^2 \right\} \right]^{1/2}$$

(incorporating a term due to the mean-correction, which is of smaller order). Using (6.14),

$$(7.37) \qquad \sum_t \phi_k(\varepsilon_t) \delta_{\ell t} = \sum_t \phi_k(\varepsilon_t) \phi_\ell'(\varepsilon_t) d_t + \frac{1}{2} \sum_t \phi_k(\varepsilon_t) \phi_\ell''(\bar{\varepsilon}_t) d_t^2.$$

We have

$$E \left\| \sum_t \{ \phi_k(\varepsilon_t) \phi_\ell'(\varepsilon_t) - E \phi_k(\varepsilon_0) \phi_\ell'(\varepsilon_0) \} d_{1t} \right\|^2 \leq CE \{ \phi_k(\varepsilon_0) \phi_\ell'(\varepsilon_0) \}^2 \sum_t E d_{1t}^2$$

$$\leq C\ell^2 \mu_{2\kappa(k+\ell+K)} \log n.$$

Replacing $d_{1t}$ by $d_t - d_{1t}$ gives no greater bound, by virtue of (6.15) and (6.17). On the other hand,

$$\{ E \phi_k(\varepsilon_0) \phi_\ell'(\varepsilon_0) \} \sum_t d_t = O_p(\ell \mu_{2\kappa k}^{1/2} \mu_{2\kappa(\ell+K)}^{1/2} n^{1/2})$$

because $\sum_t E_t = 0$ implies $\sum_t d_t = \sum_t \varepsilon_t$. Next

$$\left| \sum_t \phi_k(\varepsilon_t) \phi_\ell''(\bar{\varepsilon}_t) d_t^2 \right| \leq C^{\kappa\ell+1} \ell^2 \sum_t |\phi_k(\varepsilon_t)| (1 + |\varepsilon_t|^{\kappa(\ell+K)} + |d_t|^{\kappa(\ell+K)}) d_t^2.$$

Proceeding as in Section 6, this is $O_p((C\ell)^{2\kappa\ell+2} \mu_{\kappa k} \mu_{r_\ell} \log n)$, where $r_\ell$ is the smallest even integer exceeding $\kappa(\ell + K) + 2$. It follows that

$$\sum_{k,\ell=1}^L \left( \sum_t \phi_k(\varepsilon_t) \delta_{\ell t} \right)^2 = O_p(\rho_{2\kappa L}^2 (L^2 n + (CL)^{4\kappa L+4} (\log n)^2)).$$

Also

$$\left\{ \sum_{k,\ell=1}^L \sum_t \left( \sum_t \delta_{k\ell} \delta_{\ell t} \right)^2 \right\}^{1/2} \leq \sum_{\ell=1}^L \sum_t \delta_{\ell t}^2 \leq \sum_{\ell=1}^L \sum_t \phi_\ell'(\bar{\varepsilon}_t)^2 d_t^2,$$

and by proceeding as before this is $O_p((CL)^{4\kappa L+2} \rho_{2\kappa L} \log n)$. The proof is completed by application of Lemmas 8 and 10. $\square$

**Acknowledgments.** I thank Fabrizio Iacone for carrying out the computations reported in Section 4, and two referees for numerous comments which have led to an improved presentation.

Department of Economics
London School of Economics
Houghton Street
London WC2A 2AE
United Kingdom
e-mail: p.m.robinson@lse.ac.uk